\documentclass[11pt]{article}

\usepackage[german,english]{babel}
\usepackage[cp850]{inputenc}
\usepackage{psfrag}
\usepackage{subfigure}
\usepackage{latexsym,graphicx}
\usepackage{amsmath}
\usepackage{multirow,multicol}
\usepackage{listings}
\usepackage{color}
\font\tenmath=msbm10 scaled 1200
\font\sevenmath=msbm7 scaled 1200
\font\fivemath=msbm5 scaled 1200

\newfam\mathfam \textfont\mathfam=\tenmath
\scriptfont\mathfam=\sevenmath \scriptscriptfont\mathfam=\fivemath
\def\math{\fam\mathfam}
\def\R{{\math R}}

\def\E{{\math E}}

\def\P{{\math P}}

\def \^#1{\if#1i{\accent"5E\i}\else{\accent"5E#1}\fi}

\def \ind {1 \mkern -5mu \hbox{I}}

\def \a{\alpha}

\def \g{\gamma}

\def \ds{\displaystyle}
\def \cqfd{\quad_\diamondsuit}
\def \ms{\medskip}
\def \ss{\smallskip}
\def \bs{\bigskip}
\def \ni{\noindent}

\newtheorem{Thm}{Theorem}

\newtheorem{Lem}{Lemma}
\newtheorem{Pro}{Proposition}

\newcommand{\delimleft}[2]{\ifcase #1\or
     #2\or %
     \bigl#2\or %
     \Bigl#2\or %
     \biggl#2\or %
     \Biggl#2\or %
     \left#2\fi}
\newcommand{\delimright}[2]{\ifcase #1\or
     #2\or %
     \bigr#2\or %
     \Bigr#2\or %
     \biggr#2\or %
     \Biggr#2\or %
     \right#2\fi}
\newcommand{\pa}[2][6]{
     \delimleft{#1}{(} #2 \delimright{#1}{)}}
\newcommand{\br}[2][6]{
     \delimleft{#1}{[} #2 \delimright{#1}{]}}
\newcommand{\ac}[2][6]{
     \delimleft{#1}{\{} #2 \delimright{#1}{\}}}
\newcommand{\prob}[2][6]{
    \P\!\br[#1]{#2}}

\newcommand{\probc}[3][6]{
    \P\!\br[#1]{#2%
    \delimleft{#1}{.}%
    \vphantom{#2}\vphantom{#3}%
    \delimright{#1}{\vert} #3}}
\newcommand{\esp}[2][6]{
    \E\!{}\br[#1]{#2}}

\newcommand{\espc}[3][6]{
    \E\!{}\br[#1]{#2%
    \delimleft{#1}{.} \vphantom{#2}\vphantom{#3}%
    \delimright{#1}{\vert} #3}}

\newcommand{\psca}[2][1]{
     \delimleft{#1}{\langle} #2%
     \delimright{#1}{\rangle}}
\newcommand{\pscaLTp}[3][1]{
     \delimleft{#1}{\langle} #3%
     {\delimright{#1}{\rangle}}_{\scriptscriptstyle{\! L^2_{T,#2}}}}
\newcommand{\abs}[2][1]{
     \delimleft{#1}{\lvert} #2%
     \delimright{#1}{\rvert}}
\newcommand{\norm}[2][6]{
     \delimleft{#1}{\lVert} #2%
     \delimright{#1}{\rVert}}
\newcommand{\trnorm}[2][6]{
\delimleft{#1}{\lvert}\hspace*{-0.1em}\delimleft{#1}{\lVert} #2%
     \delimright{#1}{\rVert}\hspace*{-0.1em}\delimright{#1}{\rvert}}
\newcommand{\normsup}[2][6]{
     \delimleft{#1}{\lVert} #2%
     \delimright{#1}{\rVert}_{\scriptscriptstyle{\infty}}}
\newcommand{\normLp}[3][6]{
     \delimleft{#1}{\lVert} #3%
     \delimright{#1}{\rVert}_{\scriptscriptstyle{\!#2}}}
\newcommand{\normLTp}[3][6]{
    \normLp[#1]{L^2_{T,#2}}{#3}}
\newcommand{\indic}[2][6]{
    \boldsymbol{1}_{\ac[#1]{#2}}}

\DeclareMathOperator{\Tr}{Tr\!}
\DeclareMathOperator{\DD}{D\!}
\DeclareMathOperator{\dd}{d\!}
\DeclareMathOperator{\argmin}{Argmin}
\DeclareMathOperator{\NIG}{NIG}
\newcommand{\Dp}[1][2]{\DD^{\,#1}\!}

\newcommand{\Thetath}{\Theta^{(\theta)}}
\newcommand{\Thetamth}{\Theta^{(-\theta)}}
\newcommand{\Xthet}{X^{(\theta)}}
\newcommand{\Xthets}{X^{(\theta),s}}
\newcommand{\Xthett}{X^{(\theta),t}}
\newcommand{\Xmthet}{X^{(-\theta)}}
\newcommand{\Xmthets}{X^{(-\theta),s}}
\newcommand{\Xmthett}{X^{(-\theta),t}}
\newcommand{\bX}{\bar{X}}
\newcommand{\tk}{{t_k}}
\newcommand{\tkp}{{t_{k+1}}}

\author{{\sc Vincent Lemaire and Gilles Pag\`es}
\thanks{Laboratoire de Probabilit\'es et Mod\`eles al\'eatoires,
UMR~7599, Universit\'e Pierre et Marie Curie, case 188, 4, pl. Jussieu, F-75252
Paris Cedex 5, France. E-mail:{\tt vincent.lemaire@upmc.fr}, {\tt gilles.pages@upmc.fr}}
}
\date{\today}
\title{Unconstrained Recursive Importance Sampling}

\begin{document}
\lstset{language=C, 
    basicstyle=\ttfamily\scriptsize,
    keywordstyle=\bfseries,
    commentstyle=\slshape,
    frameround = fttt,
    emptylines={1},
    showlines= true,
    showspaces=false,
    showstringspaces=false}
\maketitle

\begin{abstract}
We propose an \emph{unconstrained} stochastic approximation method of finding the optimal measure change (in an \emph{a priori} parametric family) for Monte Carlo simulations. We consider different parametric families based on the Girsanov theorem and the Esscher transform (or exponential-tilting). In a multidimensional Gaussian framework, Arouna uses a projected Robbins-Monro procedure to select the parameter minimizing the variance (see \cite{ARO}). In our approach, the parameter (scalar or process) is selected by a classical Robbins-Monro procedure without projection or truncation.  To obtain this unconstrained algorithm we intensively use the regularity of the density of the law without assume smoothness of the payoff. We prove the convergence for a large class of multidimensional distributions and diffusion processes.

We illustrate the effectiveness of our algorithm via pricing a Basket payoff under a multidimensional NIG distribution, and pricing a barrier options in different markets.
\end{abstract}


\ni {\em 2000 Mathematics Subject Classification: 65C05, 65B99, 60H35.}

\noindent
\section{Introduction}
The basic problem in Numerical Probability is to {\em optimize} some way or another the computation by a  Monte Carlo simulation of a real quantity $m$ known by a probabilistic representation   
\begin{equation*}
    m = \esp{F(X)}
\end{equation*}
where $X:(\Omega,{\cal A}, \P)\to (E,|\,.\,|_E)$ is a random vector having values in a Banach space $E$ and $F:E\to \R$ is a Borel function (and $F(X)$ is square integrable). The space
$E$ is $\R^d$ but can also be a functional space of paths of a process $X=(X_t)_{t\in [0,T]}$. However, in this introduction section, we will first focus on the finite dimensional case $E=\R^d$. 

Assume that $X$ has an absolutely continuous distribution $\P_{_X}(dx)=p(x)\lambda_d(dx)$ ($\lambda_d$ denotes the Lebesgue measure on $(\R^d,{\cal B}or(\R^d))$) and that $F\!\in L^2(\P_{_X})$ with $\P(F(X)\neq 0)>0$ (otherwise the expectation is clearly $0$ and the problem is meaningless). Furthermore we assume that the probability density $p$ is {\em everywhere positive} on $\R^d$.  

The paradigm of importance sampling applied to  a parametrized family of distributions is the following:  consider the family of  absolutely continuous probability distributions $\pi_\theta(dx):=p_\theta(x) dx$, $\theta \!\in \Theta$, such that $p_\theta(x) > 0$, $\lambda_d(dx)$-$a.e.$ . One may assume without loss of generality that $\Theta$ is an open non empty connected subset of $\R^q$ containing  $0$ so that $p_0=p$. In fact we will assume throughout the paper that $\Theta=\R^q$. Then for any $\R^d$-valued random variable $\Xthet$ with distribution $\pi_\theta$, we have 
\begin{equation}\label{CamMar}
    \esp{F(X)} = \esp{F(\Xthet)\frac{p(\Xthet)}{p_\theta(\Xthet)}}.
\end{equation}

Among all these random variables having the same expectation $m=\esp{F(X)}$, the one with the lowest variance is the one with the lowest quadratic norm: minimizing the variance amounts to finding the parameter $\theta^* $ solution (if any) to the following minimization problem
\begin{equation*}
    \min_{\theta\!\in \R^q} V(\theta)
\end{equation*}
where, for every  $\theta\!\in \R^q$, 
\begin{equation}\label{V=E}
    V(\theta) := \esp{F^2(\Xthet) \frac{p^2(\Xthet)}{p_\theta^2(\Xthet)}} 
    = \esp{F^2(X) \frac{p(X)}{p_\theta(X)}} \le +\infty.
\end{equation}
A typical situation is importance sampling by  mean translation in a finite dimensional Gaussian framework \emph{i.e.} 
\begin{equation*}
    \Xthet = X+\theta,\quad p(x) = \frac{e^{-\frac{\abs{x}^2}{2}}}{(2\pi)^{\frac d2}}, \quad p_\theta(x)=p(x-\theta) \quad \mbox{ and }\quad V(\theta)= e^{-\abs{\theta}^2}\esp{F^2(X) e^{-2\psca{\theta,X}}}. 
\end{equation*}
Then the second equality in~(\ref{V=E})  is simply the Cameron-Martin formula.  This specific framework is very important for applications, especially in Finance,  and was the starting point of the new interest for recursive importance sampling procedures, mainly initiated by Arouna in~\cite{ARO} (see further on). In fact, as long as variance reduction  is concerned, one can consider a more general framework  without extra effort. As a matter of fact, if the distributions $p_{\theta}$ satisfy 
\begin{align}\label{H1} \tag{${\cal H}_1$}
        \begin{cases}
    (i) &     \forall \,x\! \in \R^d,\quad \theta \mapsto p_{\theta} (x) \text{ is $\log$-concave} \\ 
        (ii) & \forall\, x\! \in \R^d,\quad \lim_{\abs{\theta} \rightarrow +\infty} p_{\theta}(x) = 0,
        \qquad \text{or}\qquad
        \forall x \in \R^d,\quad \lim_{\abs{\theta} \rightarrow +\infty} \frac{p_{\theta}(x)}{p^2_{\theta/2}(x)} = 0,
        \end{cases}
\end{align}
and $F$ satisfies $\esp{F^2(X)\frac{p(X)}{p_\theta(X)}} < +\infty$ for every $\theta\!\in \R^q$, then (see Proposition~\ref{prop_intro} below), the function $V$ is finite, convex, goes to infinity at infinity. As a consequence $\argmin V=\{\nabla V=0\}$ {\em  is non empty}. Assumption $(ii)$ can be localized by by considering that one  the two conditions holds only on a Borel set $C$ of $\R^d$ such that $\P_{_X}(C\cap\{F\neq 0\})>0$.   If $\theta\mapsto p_\theta(x)$ is {\em strictly} $\log$-concave  for every $x$ in a Borel set $B$ such that $\P{_X}\!\br{B\cap\{F\neq 0\}}>0$, then $V$ is strictly convex and $\argmin V=\{\nabla V=0\}$ is reduced to a single $\theta^*\!\in \R^q$.  These results  follow from the second representation of $V$ as an expectation in~(\ref{V=E}) which  is obtained  by a second change of probability (the reverse one).  For notational convenience we will  temporarily assume that $\argmin V=\{\theta^*\}$ in this introduction section, although our main result needs no such restriction.

A classical procedure to approximate $\theta^*$ is the so-called Robbins-Monro algorithm. This is a recursive stochastic algorithm (see~(\ref{RM}) below) which can be seen as a stochastic counterpart of deterministic recursive zero search procedures  like  the Newton-Raphson one. It can be formally implemented provided the gradient of the (convex) target function $V$ admits a representation as an expectation. Since we have no {\em a priori} knowledge about the regularity of $F$ (\footnote{When $F$ is smooth enough alternative approaches have been developed based on some large deviation estimates which provide a good approximation of $\theta^*$ by deterministic optimization methods (see~\cite{GLHESH}).}) and do not wish to have any, we are naturally lead to {\em formally} differentiate  the second representation  of $V$ in~(\ref{V=E}) to obtain a representation of $\nabla V$ as
\begin{equation} \label{grad1}
    \nabla V(\theta) = \esp{F^2(X)\frac{p(X)}{p_\theta(X)} \frac{\nabla_{\!\theta} \,p_{\theta}(X)}{p_{\theta}(X)}}.
\end{equation}

Then, if we consider the function $\bar{H}_V(\theta, x)$ such that $\nabla V(\theta) = \E\left( \bar{H}_V(\theta, X) \right)$ naturally defined by~(\ref{grad1}), the derived Robbins-Monro procedure writes 
\begin{equation}\label{RM} \tag{AlgoRM}
    \theta_{n+1} = \theta_n - \gamma_{n+1} \bar{H}_V(\theta_n, X_{n+1}),
\end{equation}
with $(\gamma_n)_{n \ge 0}$ a {\em step} sequence   decreasing  to 0 (at an appropriate rate), $(X_n)_{n \ge 0}$ a sequence of i.i.d. random variables with distribution $p(x) \lambda_d(d x)$. To establish the convergence of a Robbins-Monro procedure to $\theta^* = \argmin V$ requires  seemingly not so stringent assumptions. We mean by that: not so different from those needed in a deterministic framework. However, one of them turns out to be quite restrictive for our purpose: the sub-linear growth assumption in quadratic mean
\begin{equation} \label{NEC} \tag{NEC}
    \forall \theta\!\in \R^d, \qquad \normLp{2}{\bar{H}_{V}(\theta, X)} \le C (1 + \abs{\theta}).
\end{equation}
which  is  the stochastic counterpart  of the  classical non-explosion condition needed  in a deterministic framework. In practice, this condition is almost never satisfied in our framework due  to the behaviour of the term  
$\frac{p(x)}{p_{\theta}(x)}$ as $\theta$ goes to infinity. 

\ms
The origin of recursive importance sampling as briefly described above goes back to Kushner and has  recently been  brought back to light in a Gaussian framework by Arouna in~\cite{ARO}. However, as confirmed by the numerical experiments carried out  by several authors (\cite{ARO, KAW, LEL}), the regular Robbins-Monro procedure~(\ref{RM})  does suffer from a structural instability coming from  the violation of~(\ref{NEC}). This phenomenon is quite similar to the behaviour of the explicit discretization schemes of an $ODE$ $\equiv \; \dot x=h(x)$ when $h$ has a super-linear growth at infinity. Furthermore,  in a probabilistic framework no ``implicit scheme"  can be devised  in general. Then the only way out {\em mutatis mutandis}  is to kill the procedure when it comes close to explosion and to restart it with a smaller step sequence. Formally, this can be described as some  repeated projections or truncations when the algorithm leaves a slowly growing compact set waiting for    stabilization  which is shown  to occur $a.s.$. 
Then, the algorithm behaves like a regular Robbins-Monro procedure. This is the so-called  ``Projection \`a la Chen" avatar of the Robbins-Monro algorithm, introduced by Chen in~\cite{CHE, CHEetal} and then investigated by several authors  (see $e.g.$ \cite{ANMOPR,LEL}) Formally, repeated projections ``\`a la Chen"  can be written as follows:
\begin{equation*} \tag{AlgoP}
     \theta_{n+1} = \Pi_{K_{\sigma(n)}} \ac{\theta_n - 
\bar{H}_V(\theta_n, X_{n+1})}
\end{equation*}
where $\Pi_{K_{\sigma(n)}}$ denotes the projection on the convex compact $K_{\sigma(n)}$ ($K_p$ is increasing to $\R^d$ as $p\to \infty$).  In~\cite{LEL} is established a  a Central Limit Theorem   for this version  of the recursive variance reduction procedure. Some extensions to non Gaussian framework  have been carried out by Arouna in his PhD thesis (with some applications to reliability) and more recently to the marginal distributions of a L\'evy processes by Kawai in~\cite{KAW}. 

However, convergence occurs for this procedure after a long ``stabilization phase" \dots provided that the sequence of compact sets have been specified in an appropriate way. This specification turns out to be a rather sensitive phase of the ``tuning" of the algorithm to be combined with that of the step sequence.  

\ms
In this paper, we show  that as soon as the growth  of $F$ at infinity  can be explicitly controlled, it is always possible to design a regular Robbins-Monro algorithm which  $a.s.$ converges  to a variance minimizer $\theta^*$ with no risk  of explosion (and subsequently no need of repeated projections).    

To this end the key is to introduce a {\em third} change of probability in order to control the term $\frac{p(x)}{p_{\theta}(x)}$. In a Gaussian framework this amounts to switching   the  parameter $\theta$ from the density $p$  to the function $F$ by a third   mean translation. This of course corresponds to  a new function $\bar{H}_V$ but can also be interpreted {\em a posteriori} as a way to introduce an {\em adaptive} step sequence (in the spirit of~\cite{LEM}). 

In terms of formal importance sampling, we introduce a new positive density $q_\theta$ (everywhere positive on $\{p>0\}$) so that  the gradient writes 
\begin{equation} \label{grad2}
    \nabla V(\theta) = \esp{F^2(\widetilde \Xthet) \frac{p^2(\widetilde \Xthet)}{p_\theta(\widetilde \Xthet) q_\theta(\widetilde \Xthet)} \frac{\nabla p_{\theta}(\widetilde \Xthet)}{p_{\theta}(\widetilde \Xthet)}} = \esp{\widetilde H_V(\theta,\widetilde \Xthet)},
\end{equation}
where $\widetilde \Xthet \sim q_{\theta}(x) d x$. The ``weight" $\displaystyle  \frac{p^2(\widetilde \Xthet)}{p_\theta(\widetilde \Xthet) q_\theta(\widetilde \Xthet)} \frac{\nabla p_{\theta}(\widetilde \Xthet)}{p_{\theta}(\widetilde \Xthet)}$  may seem  complicated but the r\^ole of the density $q_\theta$ is to control the critical term $\frac{p^2(x)}{p_\theta(x) q_\theta(x)}$ by a (deterministic) quantity only depending on $\theta$. Then we can replace $\widetilde H_V$ by a function $H(\theta,x)= \delta(\theta)\, \widetilde  H_V(\theta,x)$ in the above Robbins-Monro procedure~(\ref{RM}) where   $\delta $ is a positive function used to control the behaviour of $\widetilde  H_V(\theta,x)$ for large values of $x$ (note that $\ac[2]{\esp[2]{H(.,\widetilde \Xthet)} = 0} = \ac{\nabla V = 0}$). 

\medskip
We will first illustrate this paradigm in a finite dimensional setting with parametrized importance sampling procedures: the mean translation and the Esscher transform which coincide for Gaussian vectors on which a  special emphasis will be put.  Both cases correspond to a specific choice of $q_\theta$ which significantly simplifies the expression of the weight.  
 
\medskip As a second step, we will deal with an infinite dimensional setting (path-dependent diffusion like processes) where we will rely on the Girsanov transform to play the  role of mean translator. To be more precise, we want now to compute $\esp{F(X)}$ where $X$ is a path-dependent diffusion process and $F$ is a functional defined on the space ${\cal C}([0,T],\R^d)$ of continuous functions defined on $[0,T]$. We consider a $d$-dimensional It\^o process $X=(X_t)_{t\in [0,T]}$ solution of the path-dependent SDE
\begin{equation*} \tag{$E_{b,\sigma,W}$}
    \dd X_t= b(t,X^t) \dd t + \sigma(t,X^t) \dd W_t,\quad X_0=x\!\in \R^d,
\end{equation*}
where $W=(W_t)_{t\in [0,T]}$ is a $q$-dimensional standard Brownian motion, $X^t:=(X_{t\wedge s})_{s\in [0,T]}$ is the stopped process at time $t$, $b:[0,T]\times {\cal C}([0,T],\R^d)\to \R^d$ and $\sigma:[0,T]\times  {\cal C}([0,T],\R^d)\to {\cal M}(d,q)$ are Lipschitz with respect to the  $\normsup{\;.\;}$ on the space  $ {\cal C}([0,T],\R^d)$ and continuous in $(t,x)\!\in [0,T]\times {\cal C}([0,T],\R^d)$ (see~\cite{ROWI} for more details about these path-dependent SDE's).

Let $\varphi$ be a fixed borel bounded functional on ${\cal C}([0,T],\R^d)$ with values in ${\cal M}(q, p)$  (where $p \ge 1$ is a free integral parameter). Then a Girsanov transform yields that for every $\theta\!\in L^2_{T,p}:=L^2([0,T],\R^p)$, 
\begin{equation*}
    \esp{F(X)} = \esp{F(\Xthet) e^{-\int_0^T\!\psca{\varphi(\Xthets) \theta(s), \dd W_s} - \frac{1}{2} \normLTp{q}{\varphi(\Xthet{}^{,.})\theta}^2}}
\end{equation*}
where $\Xthet$ is the solution to $(E_{b+\sigma\varphi\theta,\sigma})$. The functional to be minimized is now 
$$
V(\theta) = \esp{F(\Xthet)^2 e^{-2\int_0^T\!\psca{\varphi(\Xthets) \theta(s), \dd W_s} -  \normLTp{q}{\varphi(\Xthet{}^{,.})\theta}^2}},\qquad \theta\!\in L^2_{T, p}.
 $$
 In practice we will only minimize $V$ over a {\em finite dimensional subspace} of $E={\rm span}\{e_1,\ldots,e_m\}\subset L^2_{T, p}$. 

The paper is organized as follows. Section \ref{dimfinie} is devoted to the finite dimensional setting where we recall the main tool including a slight extension of the Robbins-Monro theorem in the Subsection \ref{argmin-target} and the gaussian case investigated in \cite{ARO} is revisited to emphasize the new aspects of our algorithm in the Subsection \ref{gaussian_case}.

In Section 2 we successively investigate the translation for log-concave distributions probability and the Esscher transform. In Section \ref{diffusion} we introduce a functional version of our algorithm based on the Girsanov theorem to deal the SDE.
In Section 4 we provide some comments on the practical implementation and in Section 
\ref{numerical} some numerical experiments are carried out on some option pricing problems.

\ss
\ni {\sc Notations:} $\bullet$ We will denote by $S>0$ the fact that a symmetric matrix $S$  is positive definite. 
\ni $|\,.\,|$ will denote the canonical Euclidean norm on $\R^m$ and $\langle\,.,\,.\,\rangle$ will denote the canonical inner product.

\ni $\bullet$ The real constant $C>0$ denotes a positive real constant that may vary from line to line.

\ni $\bullet$ $\normLTp{p}{f}:= \left(\int_0^T f_1^2(t)+\cdots+f_p^2(t)dt \right)^{\frac 12}$ if $f=(f_1,\ldots,f_p)$ is an $\R^p$-valued (class of) Borel function(s).

\section{The finite-dimensional setting} \label{dimfinie}
\subsection{$\boldsymbol \argmin \,\boldsymbol V$ as a target} \label{argmin-target}
\begin{Pro} \label{prop_intro}
Suppose \eqref{H1} holds. \\ Then the function $V$ defined by~(\ref{V=E}) is  convex and $\lim_{\abs{\theta} \rightarrow +\infty} V(\theta) = +\infty$. As a consequence
\[
    \argmin V= \{\nabla V=0\}\neq \emptyset.
\] 
\end{Pro}
\noindent {\bf Proof.}
By the change of probability $\frac{d \pi_{\theta}}{d \lambda_d}$ we have $V(\theta) = \esp{F^2(X)\frac{p(X)}{p_\theta(X)}}$. Let $x$ fixed in $\R^d$. The function $(\theta \mapsto \log p_\theta(x))$ is concave, hence $\log (1/p_\theta(x))=-\log p_\theta(x)$ is convex so that, owing to the Young Inequality, the function $\frac{1}{p_\theta(x)}$ is convex since it is non-negative.

To prove that $V$ tends to infinity as $\abs{\theta}$ goes to infinity, we consider two cases:
\begin{itemize}
    \item[--] If $\ds \lim_{\abs{\theta} \rightarrow +\infty} p_\theta(x) = 0$ for every $x \in \R^d$, the result is trivial by Fatou's Lemma.
    \item[--] If $\ds \lim_{\abs{\theta} \rightarrow +\infty} \frac{p_\theta(x)}{p^2_{\theta/2}(x)} = 0$ for every $x \in \R^d$, we apply the reverse H\"older inequality with conjugate exponents $(\frac{1}{3},-\frac{1}{2})$ to obtain 
\begin{align*}
    V(\theta) & \ge \esp[5]{F^{2/3}(X) \pa[4]{\frac{p^2_{\theta/2}(X)}{p(X)p_\theta(X)}}^{\frac{1}{3}}}^3 \esp{\pa{\frac{p(X)}{p_{\theta/2}(X)}}^{-1}}^{-2}, \\
    & \ge \esp[5]{F^{2/3}(X) \pa[4]{\frac{p^2_{\theta/2}(X)}{p(X)p_\theta(X)}}^{\frac{1}{3}}}^3,
\end{align*}
($p$ and  $p_\theta$ are probability   density functions).
One concludes again by Fatou's Lemma. $\cqfd$
\end{itemize}

\ms
The set $\argmin V$, or to be precise, the random vectors taking values in $\argmin V$ will the target(s) of our new algorithm. If $V$ is strictly convex, $e.g.$ if $$\prob{X\!\in\{p_.(x) \mbox{ strictly $\log$-concave and } F(x)\neq 0\}} > 0,$$ then $\argmin V=\{\theta^*\}$. Nevertheless this will not be necessary owing to the combination of the two results that follow. 


\begin{Lem}\label{Un}Let $U:\R^d\to \R_+$ be a convex differentiable function, then
\begin{equation*}
    \forall\, \theta,\, \theta'\!\in \R^d,\qquad \psca{\nabla U(\theta)-\nabla U(\theta'), \theta-\theta'} \ge 0.
\end{equation*}
Furthermore, if $\argmin U$ is nonempty, it is a convex closed set (which coincide with $\{\nabla U=0\}$) and 
\begin{equation*}
    \forall\, \theta\!\in \R^d\setminus \argmin U, \;\forall\, \theta^*\!\in \argmin U, \qquad \psca{\nabla U(\theta), \theta-\theta^*}> 0.
\end{equation*}
\end{Lem}

A sufficient (but in no case necessary) condition for a nonnegative convex function $U$ to attain a minimum is that $\lim_{|x|\to \infty} U(x)=+\infty$. 
  
\medskip
Now we pass to the statement of the convergence theorem on which we will rely throughout the paper. It is a slight variant of the regular Robbins-Monro procedure whose proof is rejected in an annex.

\begin{Thm}\label{ThmRZ} (Extended Robbins-Monro Theorem)
Let $H:\R^q\times\R^d\to \R^d$ a Borel function and $X$ an $\R^d$-valued random vector such that $\esp{\abs{H(\theta,X)}}<+\infty$ for every $\theta\!\in \R^d$. Then set
\begin{equation*}
    \forall\, \theta\!\in \R^d,\qquad h(\theta)= \esp{H(\theta,X)}.
\end{equation*}
Suppose that the function $h$ is continuous and that ${\cal T}^*:=\{h=0\}$ satisfies
\begin{equation}\label{RMmeanreverting}
    \forall\, \theta\!\in \R^d\setminus {\cal T}^* ,\; \forall\, \theta^*\!\in{\cal T}^*,  \qquad \psca{\theta-\theta^*,h(\theta)} > 0.
\end{equation}
Let $\g=(\g_n)_{n\ge 1}$ be a sequence of gain parameters satisfying
\begin{equation}\label{StepCond}
  \sum_{n\ge 1} \g_n=+\infty \qquad \mbox{ and }\qquad \sum_{n\ge 1} \g^2_n<+\infty.
\end{equation}
Suppose that
\begin{equation}\label{LinGrowth} \tag{NEC}
    \forall\, \theta\!\in \R^d,\qquad \esp{\abs{H(\theta,X)}^2} \le C(1+\abs{\theta}^2)
\end{equation}
(which implies $\abs{h(\theta)}^2\le C(1+\abs{\theta}^2)$).

Let $(X_n)_{n\ge 1}$ be an i.i.d. sequence of random vectors having the distribution of $X$, a random vector  $\theta_0$, independent of $(X_n)_{n\ge 1}$   satisfying $\esp{|\theta_0|^2} <+\infty$, all defined on the same probability space $(\Omega,{\cal A}, \P)$. Then, the recursive procedure defined by
\begin{equation}\label{Algo}
    \theta_n = \theta_{n-1}-\g_{n+1} H(\theta_n,X_{n+1}),\qquad n \ge 1
\end{equation}
satisfies:
\begin{equation}
    \exists \, \theta_{_\infty}:(\Omega,{\cal A})\to {\cal T}^*,\; \theta_{_\infty}\!\in L^2(\P),\quad \mbox{such that} \quad \theta_n  \stackrel{a.s.}{\longrightarrow} \theta_{_\infty}.\hskip 1,5 cm 
\end{equation}
The convergence also holds in  $L^p(\P)$, $p\!\in (0,2)$.
\end{Thm}

The proof is postponed to the Appendix at the end of the paper. The natural way to apply this theorem for our purpose is the following: 
\begin{itemize} 
    \item[--] {\sc Step~1}:   we will show that the convex function $V$ in~(\ref{V=E}) is differentiable with a gradient  $\nabla V$   having a representation as an expectation formally given $\nabla V (\theta)= \esp{\nabla_{\theta} v(\theta,X)}$.
    \item[--] {\sc Step~2}:  then set $H(\theta,x):= \rho(\theta)\nabla_{\theta} v(\theta,x)$ where $\rho$ is a (strictly) positive function on $\R^q$. As a matter of fact, with the notations of the above theorem
\begin{equation*}
    \psca{\theta-\theta^*, h(\theta)} = \delta(\theta) \psca{\theta-\theta^*, \nabla V(\theta)},
\end{equation*}
so that ${\cal T}^*= \argmin V$ and \eqref{RMmeanreverting} is satisfied (set $U:=V$ in Lemma~\ref{Un}). 
    \item[--] {\sc Step~3}:  Specify in an appropriate way the function $\delta$ so that the linear quadratic growth assumption is satisfied. This is the sensitive point that will lead us   to modify   the structure more deeply  by finding a new representation of $\nabla V$ as an expectation not directly based on the  local gradient $\nabla_{\theta} v(\theta,x)$.
\end{itemize}

\subsection{A first illustration: the Gaussian case revisited} \label{gaussian_case}
The Gaussian is the framework of~\cite{ARO}. It is also a kind of introduction to the infinite dimensional diffusion setting investigated in Section~\ref{diffusion}. In the Gaussian case, the natural importance sampling density is the translation of the gaussian density: $p_{\theta}(x) = p(x-\theta)$ for $\theta \in \R^d$ ($i.e.$ $q=d$). We have
\begin{equation*}
     p_\theta(x) = e^{-\frac{\abs{\theta}^2}{2} + \psca{\theta, X}} p(x).
\end{equation*}

The assumption \eqref{H1} is clearly satisfied by the Gaussian density, and we assume that $F$ satisfies $\esp{   F^2(X) e^{- \psca{\theta, X}}} < +\infty$ so that $V$ is well defined.

In~\cite{ARO}, Arouna considers the function $\bar{H}_V(\theta, x)$ defined by
\begin{equation*}
     \bar{H}_V(\theta, x) = F^2(x) e^{\frac{\abs{\theta}^2}{2} - \psca{\theta, x}}(\theta-x).
\end{equation*}
It is clear that the condition \eqref{NEC} is not satisfied  even if we simplify this function by $e^{\abs{\theta}^2/2}$ (which does not modify the problem).

\medskip 
\ni {\em A first approach:} When $F(X)$ have  finite moments of any order, a naive way to control  directly $\normLp{2}{\bar{H}_V(\theta_n, X_{n+1})}$ by an explicit  deterministic function of $\theta$ (in order to rescale it) is to proceed as follows: one derives from H\"older Inequality that for every couple $(r,s)$, $r,\,s>1$ of conjugate exponents
\begin{equation*}
     \normLp{2}{\bar{H}_V(\theta, X)} \le 
        \normLp{2r}{(\theta-X)F^2(X)} \normLp{2s}{e^{-\psca{\theta, X}}} e^{\frac{\abs{\theta}^2}{2}}.
\end{equation*}
Setting $r=1+\frac{1}{\varepsilon}$ and $s=1+\varepsilon$,  yields
\begin{align*}
     \normLp{2}{\bar{H}_{V}(\theta, X)}
     & \le\pa{\normLp{2(1+\frac{1}{\varepsilon})}{X F^2(X)} + 
        \normLp{2(1+\frac{1}{\varepsilon})}{F^2(X)}\abs{\theta}} 
        e^{\pa{\frac{3}{2}+\varepsilon}\abs{\theta}^2}.
\end{align*}

Then, $\bar{H}_{\varepsilon}(\theta, x):=e^{-(\frac{3}{2}+\varepsilon)\abs{\theta}^2}\bar{H}_{V}(\theta, x)$ satisfies the condition \eqref{NEC}  and theoretically the standard Robbins-Monro algorithm implemented with $\bar{H}_{\varepsilon}$ $a.s.$ converges and no  projection nor truncation is needed. 
Numerically, the solution is not satisfactory because the correcting factor $e^{-\pa{\frac{3}{2}+\varepsilon}\abs{\theta}^2}$ goes to zero much too fast as $\theta$ goes to infinity: if at any iteration at the beginning of the procedure $\theta_n$ is sent ``too far", then it is frozen instantly.  If $\varepsilon$ is too small it will simply not prevent explosion. The tuning of $\varepsilon$ becomes quite demanding and payoff dependent. This is in complete contradiction with our aim of {\em a self-controlled variance reducer}. A more robust approach  needs to be developed. On the other hand this kind of behaviour suggests that we are not in the right asymptotics to control $\normLp{2}{\bar{H}_V(\theta_n, X_{n+1})}$.

\smallskip Note however that when $F$ is bounded with a compact support, then one can set $\varepsilon=0$ and the above approach provides an efficient answer to our problem.

\medskip
\ni {\em A general approach:} We consider the density 
\begin{equation*}
    q_\theta(x) = \frac{p^2(x)}{p(x-\theta)} = p(x+\theta).
\end{equation*}
By \eqref{grad2}, we have
\begin{equation*}
     \nabla V(\theta) = \esp{F^2(\widetilde \Xthet) \frac{\nabla p(\widetilde \Xthet-\theta)}{p(\widetilde \Xthet-\theta)}},
\end{equation*}
with $\widetilde \Xthet \sim p(x+\theta) d x$, \emph{i.e.} $\widetilde \Xthet = X-\theta$. Since $p$ is the Gaussian density, we have $\frac{\nabla p(x)}{p(x)} = - x$. As a consequence, the function $H_V$ defined by
\begin{equation*}
     H_V(\theta, x) = F^2(x-\theta) (2\theta - x),
\end{equation*}
provides a representation $\nabla V(\theta) = \esp{H_V(\theta, X)}$ of the gradient  $\nabla V$. As soon as  $F$ is bounded, this function satisfies the condition \eqref{NEC}. Otherwise, we note that thanks to this new change of variable the parameter $\theta$ lies now {\em inside} the payoff function $F$ and that the {\em exponential term} has disappeared from the expectation.  If we have an {\em a priori} control on the function $F(x)$ as $|x|$ goes to infinity, say 
\begin{equation*}
     \exists \lambda \in \R_+, \quad \forall x \in \R^d, \quad \abs{F(x)} \le c_F e^{\lambda \abs{x}},
\end{equation*}
then we can consider the function $H_{\lambda}(\theta, x) = e^{-\lambda \abs{\theta}} H_V(\theta, x)$ which satisfies
\begin{align*} \begin{split}
     \normLp{2}{H_\lambda(\theta, X)} & \le c_F^2 
\normLp{2}{e^{2\lambda \abs{X}} (2 \theta - X)}, \\
     & \le C (1 + \abs{\theta}).
     \end{split}
\end{align*}

The resulting Robbins-Monro algorithm reads
\begin{equation*}
     \theta_{n+1} = \theta_n - \gamma_{n+1} e^{-\lambda \abs{\theta_n}} F^2(X_{n+1} - \theta_n) (2 \theta_n - X_{n+1}).
\end{equation*}

We  no longer to tune the correcting factor and one verifies on simulations that it does not suffer from freezing in general. In case of a too dissymmetric function $F$ this may still happen but a  {\em self-controlled} variant is proposed in Section~\ref{translation} below to completely get rid of this effect (which cannot be compared to an explosion).

\subsection{Translation of the mean: the general strongly unimodal case}\label{translation}
We consider importance sampling by mean translation, namely we set 
\begin{equation*}
     \forall x \in \R^d, \quad p_\theta(x) = p(x-\theta),
\end{equation*}
for $\theta \in \R^d$.

In this section we assume that 
\begin{equation*}
    p \mbox{ is $\log$-concave and } \lim_{|x|\to \infty} p(x)=0,
\end{equation*}
so that \eqref{H1} holds. Moreover, we make the following additional assumption on the probability density $p$
\begin{equation} \tag{${\cal H}^{tr}_a$} \label{H2}
\exists\,a\!\in[1,2]  \;  \mbox{ such that } \;\left\{
\begin{array}{ll}(i)& \frac{\left|\nabla p\right|}{p}(x)=O(|x|^{a-1}) 
\; \mbox{ as }\; |x|\to \infty \\ (ii)&\exists \delta>0,\; \log 
p(x)+\delta |x|^a\;\mbox{ is convex.}
\end{array}\right.
\end{equation}

First we will use~(\ref{V=E}) to differentiate $V$ since
\begin{Pro} Suppose \eqref{H1} and \eqref{H2} are satisfied and the function
$F$ satisfies
\begin{equation}\label{CondDiff}
    \forall\, \theta \in \R^d, \quad \esp{F^2(X) \frac{p(X)}{p(X-\theta)}} < +\infty \quad\text{ and }\quad \forall\, C>0,\quad \esp{F^2(X)e^{C|X|^{a-1}}} < +\infty.
\end{equation}
Then $V$ is finite and differentiable on $\R^d$ with a gradient given by
\begin{align}\label{gradV1}
    \nabla V(\theta) & = \esp{F^2(X)\frac{p(X)}{p^2(X-\theta)}\nabla p(X-\theta)}\\
    & = \esp{F^2(X-\theta)\frac{p^2(X-\theta)}{p(X)p(X-2\theta)}\frac{\nabla p(X-2\theta)}{p(X-2\theta)}}.\label{gradV2} 
\end{align}
\end{Pro}

\noindent {\bf Proof.} The formal differentiation to get~(\ref{gradV1}) from~(\ref{V=E})  is obvious. So it remains to check the domination property for $\theta$ lying inside a compact set. Let $x\!\in \R^d$ and $\theta\!\in \R^d$. The $\log$-concavity of $p$ implies that
\[
\log p(x)\le \log p(x-\theta) +\frac{\psca{\nabla p(x-\theta),\theta}}{p(x-\theta)}
\]
so that
\[
0\le \frac{p(x)}{p(x-\theta)}\le \exp{\frac{|\nabla 
p(x-\theta)|}{p(x-\theta)}|\theta|}.
\]
Using the assumption~\eqref{H2} yields, for every $\theta\!\in B(0,R)$,
\begin{align*}
F^2(X)\frac{p(X)}{p^2(X-\theta)}\nabla p(X-\theta)
&\le F^2(X)(A|X|^{a-1}+B)e^{(A|X-\theta|^{a-1}+B)|\theta|}\\
& \le C_{R}F^2(X) e^{C'|X|^{a-1}}=:g(X)\!\in L^1(\P)
\end{align*}
To derive the second expression \eqref{gradV2} for the gradient, we proceed as
follows: an elementary change of variable shows that
\begin{align*}
\nabla V(\theta)&= \int_{\R^d}F^2(x)\frac{p(x)}{p^2(x-\theta)}p(x)\nabla 
p(x-\theta)dx\\
&= \int_{\R^d}F^2(x-\theta)\frac{p^2(x-\theta)}{p^2(x-2\theta)}\nabla 
p(x-2\theta)dx\\
&= \esp{F^2(X-\theta)\frac{p^2(X-\theta)}{p(X)p(X-2\theta)}\frac{\nabla
p(X-2\theta)}{p(X-2\theta)}}.\cqfd
\end{align*}

\medskip
\ni {\bf Remark.} The second change of variable (in~\eqref{V=E}) has been processed to withdraw the parameter $\theta$ from the possible non smooth function $F$ to make possible  the differentiation of $V$ (since $p$ is smooth). The second expression~\eqref{gradV2} results form a {\em third} translation of the variable in order to plug back the parameter $\theta$ into the function $F$ which in common applications  has a known controlled growth rate at infinity.

This last statement may look strange at a first glance since $\theta$ appears in the ``weight" term of the expectation that involves the probability density $p$. However,  when $X\stackrel{d}{=}{\cal N}(0;1)$, this term can be controlled easily since it reduces to
\[
\frac{p^2(x-\theta)}{p(x)p(x-2\theta)}\frac{\nabla 
p(x-2\theta)}{p(x-2\theta)}= e^{\theta^2}(2\theta-x).
\]
The following lemma shows that, more generally in our strongly unimodal setting, if \eqref{H1} and \eqref{H2} are satisfied, this ``weight" can always be controlled  by a deterministic function of $\theta$.

\begin{Lem}\label{lemmetec} If \eqref{H2} holds, then there exists two real constants $A$, $B$ such that
\begin{equation}\label{IneqTech}
    \frac{p^2(x-\theta)}{p(x)p(x-2\theta)}\frac{|\nabla p(x-2\theta)|}{p(x-2\theta)}\le e^{2\delta |\theta|^a}\left(A |x|^{a-1} + A |\theta|^{a-1}+B\right).
\end{equation}
\end{Lem}

\noindent {\bf Proof.} Let $f$  be the convex function defined on $\R^d$ by $f(x)= \log p(x)+\delta|x|^a$. Then, for every $x,\, \theta\!\in \R^d$,
\[
\log\left(\frac{p^2(x-\theta)}{p(x)p(x-2\theta)}\right)= 
\log\left(\frac{f^2(x-\theta)}{f(x)f(x-2\theta)}\right)+ 
\delta\left(|x|^a+|x-2\theta|^a -2|x-\theta|^a\right)
\]
Note that $x-\theta= \frac 12 (x+(x-2\theta))$. Then, using the $\log$-convexity of $f$ and the elementary inequality
\[
\forall\, u,\, v\!\in \R^d,\quad |u|^a+|v|^a \le 2 
\left(\left|\frac{u+v}{2}\right|^a+\left|\frac{u-v}{2}\right|^a\right)
\]
(valid if $a\!\in (0,2]$) yields
\begin{equation} \label{control_p_trans}
     \frac{p^2(x-\theta)}{p(x)p(x-2\theta)}\le e^{2\delta |\theta|^a}.
\end{equation}
One concludes by the point $(i)$ of \eqref{H2}.  $\cqfd$

\bs
\ni {\bf Remark.} Thus the normal distribution satisfies \eqref{H2} with $a=2$ and $\delta=1/2$. Moreover, note that the last inequality in the above proof holds as an equality.

\bs Now we are in position to derive an unconstraint (extended) Robbins-Monro algorithm to minimize the function $V$, provided   the function $F$ satisfies a sub-multiplicative control property, in which $c>0$ is a real parameter and $\widetilde F$ a function from $\R^d $ to $\R_+$, such that, namely
\begin{equation} \tag{${\cal H}^{tr}_{c}$} \label{H3}
    \left\{\begin{array}{ll} 
    \forall\, x,\, y\!\in \R^d, & |F(x)|\le \widetilde F(x) \; \mbox{ and }\;
      \widetilde F(x+y)\le C(1+\widetilde F(x))^c(1+\widetilde F(y))^c\\ \\
    & \quad  \esp{|X|^{2(a-1)}\widetilde F(X)^{4c}}<+\infty .
     \end{array}\right.
\end{equation}


\noindent{\bf Remark.} 
Assumption \eqref{H3} seems almost non-parametric. However, its field of application is somewhat limited by~(\ref{H2}) for the following reason:  if  there exists a positive real number  $\eta>0$ such that  $x\mapsto \log p(x)  +\eta |x|^a$ is concave, then $p(x) \le Ce^{-\eta |x|^a}(|x|+1)$ for some real constant $C>0$; which in turn implies that the function $\widetilde F$ in \eqref{H3} needs to satisfy $\widetilde F(x)\le C'e^{\lambda |x|^b}$ for some $b\!\in (0,a)$ and some $\lambda>0$. (Then $c=c_b$ with $c_b=1$ if $b\!\in[0,  1]$ and $c_b = 2^{\frac b2}$ if $b\!\in (1,a)$, when $a>1$). 

\begin{Thm}\label{ThmMT}
Suppose $X$ and $F$ satisfy \eqref{H1}, \eqref{H2},~(\ref{CondDiff}) and  \eqref{H3} for some parameters $a\!\in (0,2]$, $b\!\in(0,a)$ and $\lambda>0$, and that the step sequence $(\gamma_n)_{n\ge 1}$ satisfies the usual decreasing step assumption
\[
\sum_{n\ge 1}\g_n=+\infty \qquad \mbox{ and } \qquad \sum_{n\ge 
1}\g^2_{n+1} <+\infty.
\]
Then the recursive procedure defined by
\begin{equation}\label{AlgoRM}
    \theta_{n+1}=\theta_n-\g_{n+1}H(\theta_n, X_{n+1}),\qquad \theta_0\!\in \R^d
\end{equation}
where   $(X_n)_{n\ge 1}$ is an i.i.d. sequence with the same 
distribution as $X$ and
\begin{equation}\label{H-a}
H(\theta,x):= 
\frac{F^2(x-\theta)}{1+\widetilde F(-\theta)^{2c}}e^{-2\delta 
|\theta|^a}\frac{p^2(x-\theta)}{p(x)p(x-2\theta)}\frac{\nabla
p(x-2\theta)}{p(x-2\theta)},
\end{equation}
  $a.s.$ converges toward an $\argmin V$-valued (square integrable) random variable $\theta^*$.
\end{Thm}

\ni {\bf Proof.} In order to apply Theorem~\ref{ThmRZ}, we have to check the following fact:

-- {\em Mean reversion}: The mean function of the procedure defined by~(\ref{AlgoRM}) reads
\begin{equation*}
    h(\theta) = \esp{H(\theta,X)} = \frac{e^{-2\delta |\theta|^a}}{1+\widetilde F(-\theta)^{2c}}\nabla V(\theta)
\end{equation*}
so that ${\cal T}^*:=\{h=0\}=\{\nabla V =0\}$ and if   $\theta^*\!\in {\cal T}^*$ and $\theta\!\in \R^d\setminus {\cal T}^*$, 
\begin{equation*}
    \psca{\theta-\theta^*, h(\theta) } =  \frac{e^{-2\delta |\theta|^a}}{1+\widetilde F(-\theta)^{2c}} \psca{\nabla V(\theta), \theta-\theta^*} >0
\end{equation*} 
for every $\theta\neq \theta^*$.

-- {\em Linear growth of $\theta \mapsto \|H(\theta,X)\|_{_2}$}: All our efforts in the design of the procedure are   motivated by this Assumption~(\ref{LinGrowth}) which prevents  explosion. This condition is clearly fulfilled by $H$ since
\begin{align*}
    \esp{\abs{H(\theta,X)}}^2 &= \frac{e^{-4\delta |\theta|^a}}{(1+\widetilde F(-\theta)^{2c})^2}\esp{F^4(X-\theta)\left(\frac{p^2(X-\theta)}{p(X)p(X-2\theta)}\frac{|\nabla p(X-2\theta)|}{p(X-2\theta)}\right)^2}, \\ 
&\le C e^{-4\delta |\theta|^a} \esp{(1+\widetilde F(X)^{2c})^2 (A(|X|^{a-1}+|\theta|^{a-1})+B)^2},
\end{align*}
where we used Assumption~\eqref{H3} in the first line and Inequality~(\ref{IneqTech}) from Lemma~\ref{lemmetec} in the second line.  One derives that there exists a real constant $C>0$ such that
\begin{equation*}
    \esp{\abs{H(\theta,X)}}^2 \le C \esp{ \widetilde F(X)^{4c}(1+|X|)^{2(a-1)}} (1+|\theta|^{2(a-1)}).
\end{equation*}
This provides the expected condition since~\eqref{H3}  holds. $\cqfd$

\subsubsection*{Examples of distributions}
\begin{itemize}
\item{The normal distribution.} Its density is given on $\R^d$ by
\[
p(x)= (2\pi)^{-\frac d2}e^{-|x|^2/2},\qquad x\!\in \R^d.
\]
so that \eqref{H1} is satisfied as well as \eqref{H2} for $a=2$, $\delta=\frac12$. Assumption \eqref{H3} is satisfied iff $(b,\lambda)\!\in(0,2)\times (0,\infty)\cup\{2\}\times (0,\frac 12)$. 
Then the function $H$ has a particularly simple form
\[
H(\theta,x) =e^{-\frac{\lambda}{2}|\theta|^b }F^2(x-\theta) (2\theta-x)
\]
\item {\em The hyper-exponential distributions} 
\[
p(x) = C_{d,a,\sigma}e^{-\frac{|x|^a}{\sigma^a}}P(x), \quad a\!\in [1,2]
\]
where $P$ is polynomial function. This wide family includes the normal distributions, the Laplace distribution, the symmetric gamma distributions, etc.

\item {\em The logistic distribution} Its density on the real line  is given by
\[
p(x) =\frac{e^x}{(e^x+1)^2}
\]
\eqref{H1} is satisfied as well as \eqref{H2} for $a=1+\eta$ ($\eta\!\in(0,1)$), $\delta>0$.  Assumption \eqref{H3} is satisfied iff $(b,\lambda)\!\in(0,1)\times (0,\infty)\cup\{1\}\times (0,1)$.  
\end{itemize}

\subsection{Exponential change of measure: the Esscher transform}\label{Esscher}
A second classical approach is to consider an exponential change of measure (or Esscher transform). This transformation has already been consider for that purpose in~\cite{KAW} to extend the procedure with repeated projections introduced in~\cite{ARO}. We denote by $\psi$ the cumulant generating function (or log--Laplace) of $X$ \emph{i.e.} $\psi(\theta) = \log \esp{e^{\psca{\theta, X}}}$. We assume that $\psi(\theta)<+\infty$ for every $\theta\!\in \R^d$ (which implies that $\psi$ is an infinitely differentiable convex function) and define
\begin{equation*}
     p_\theta(x) = e^{\psca{\theta,x} - \psi(\theta)} p(x), \quad  x \in \R^d.
\end{equation*}
Let $X^{(\theta)}$ denote any random variable with distribution $p_\theta$.

We assume that $\psi$ satisfies 
\begin{equation} \label{Hes} \tag{${\cal H}^{es}_\delta$}
    \lim_{\abs{\theta}\rightarrow +\infty} \psi(\theta)-2\psi\pa[3]{\frac{\theta}{2}} = +\infty \qquad \text{and} \qquad \exists\,\delta>0,\quad \theta \mapsto \psi(\theta) - \delta \abs{\theta}^2 \text{ is concave}.
\end{equation}

One must be aware that what follows makes sense as a variance reduction procedure only if the distribution of $\Xthet$ can be simulated at the same cost as $X$ or at least at a reasonable cost $i.e.$
\begin{equation}\label{ReXEsscher}
X^{(\theta)}= g(\theta,\xi),\quad \xi :(\Omega,{\cal A}, \P)\to {\cal X}
\end{equation}
where $ {\cal X}$ is a Borel subset of a metric space and $g:\R^d\times {\cal X}$ is an explicit Borel function.
By \eqref{V=E}, the potential $V$ to be minimized is $V(\theta) = \esp{F^2(X) e^{-\psca{\theta,X} + \psi(\theta)}}.$

\begin{Pro} Suppose $\psi$ satisfies \eqref{Hes} and $F$ satisfies 
\begin{equation} \label{hypo-F-es}
    \forall \theta \in \R^d,\quad \esp{\abs{X} F^2(X) e^{\psca{\theta,X}}} < +\infty.
\end{equation}
Then \eqref{H1} is fulfilled and the function $V$ is differentiable on $\R^d$ with a gradient given by  
\begin{align} \label{formal_diff_esscher}
    \nabla V(\theta) &= \esp{\pa{\nabla \psi(\theta) - X} F^2(X) e^{-\psca{\theta, X} + \psi(\theta)}}, \\ 
    &= \esp{\pa[2]{\nabla \psi(\theta) - \Xmthet} F^2(\Xmthet)} e^{\psi(\theta)-\psi(-\theta)}, \label{diff_esscher}
\end{align}
where $\ds \nabla \psi(\theta) = \frac{\esp{X e^{\psca{\theta,X}}}}{\esp{e^{\psca{\theta,X}}}}$.
\end{Pro}

\ni {\bf Proof.}
The function $\psi$ is clearly log-convex so that $\theta \mapsto p_\theta(x)$ is $\log$-concave for every $x\!\in \R^d$. On the other hand, by \eqref{Hes} we have $\lim \frac{p^2_{\theta/2}(x)}{p_\theta(x)} = +\infty$ for every $x\!\in \R^d$, and \eqref{H1} is fulfilled.

The formal differentiation to get \eqref{formal_diff_esscher} is obvious and is made rigorous by applying the assumption on $F$. The second expression \eqref{diff_esscher} of the gradient uses a third change of variable 
\begin{align*}
     \nabla V(\theta) & = \int_{\R^d}\! \pa{\nabla \psi(\theta) - x} F^2(x) e^{-\psca{\theta,x} + \psi(\theta)} p(x) dx, \\
     & = \int_{\R^d}\! \pa{\nabla \psi(\theta) - x} F^2(x) e^{\psi(\theta) - \psi(-\theta)} p_{-\theta}(x) dx, \\
     &= \esp{\pa[2]{\nabla \psi(\theta) - \Xmthet} F^2(\Xmthet)} e^{\psi(\theta)-\psi(-\theta)}.\cqfd
\end{align*}

\begin{Thm} We assume that $\psi$ satisfies~\eqref{Hes} and $F$ satisfies~\eqref{hypo-F-es} and
\begin{equation*}
  \forall x \in \R^d, \quad \abs{F(x)} \le C e^{\frac{\lambda}{4} \abs{x}}  \qquad \text{and}\qquad  \esp{\abs{X}^2 e^{\lambda \abs{X}}} < +\infty.
   \end{equation*}
Then the recursive procedure 
\begin{equation*}
   \left\{\begin{array}{lcl}X_{n+1}^{(\theta_n)}&=& g(\theta_n,\xi_{n+1})\\ \theta_{n+1}&=&\theta_n -\g_{n+1}H(\theta_n, X_{n+1}^{(-\theta_n)}),\quad n\ge 0, \quad  \theta_0\!\in \R^d
   \end{array}\right.
\end{equation*}
where $(\xi_n)_{n\ge 1}$ is an i.i.d. sequence with the same distribution as $\xi$  in~(\ref{ReXEsscher}) and
\begin{equation*}
    H(\theta, x) := e^{-\frac{\lambda}{2}\sqrt{d} \abs{\nabla \psi(-\theta)}} F^2(x) \pa[2]{\nabla \psi(\theta) - x}
\end{equation*}
$a.s.$ converges toward an $\argmin V$-valued (square integrable) random vector $\theta^*$ .
\end{Thm}

\ni {\bf Proof.}
We have to check the linear growth of the function $\theta \mapsto \normLp{2}{H(\theta, \Xmthet)}$ (condition \eqref{NEC}). We have
\begin{align}
    \esp{\abs{H(\theta, \Xmthet)}^2} & = e^{-\lambda\sqrt{d} \abs{\nabla \psi(-\theta)}} \esp{F^4(\Xmthet) \abs[2]{\nabla \psi(\theta) - \Xmthet}^2}, \notag \\
    & \le C e^{-\lambda\sqrt{d} \abs{\nabla \psi(-\theta)}} \esp{e^{\lambda \abs{\Xmthet}}\abs[2]{\nabla \psi(\theta) - \Xmthet}^2}, \notag \\
    & \le C e^{-\lambda \sqrt{d}\abs{\nabla \psi(-\theta)}} \pa{\abs{\nabla \psi(\theta)}^2 \esp{e^{\lambda \abs{\Xmthet}}} + \esp{\abs{\Xmthet}^2 e^{\lambda \abs{\Xmthet}}}}. \label{esscher-check-NEC}
\end{align}

First, by the following inequality
\begin{equation*}
    \forall x \in \R^d, \quad e^{\lambda \abs{x}} \le \prod_{j=1}^d \pa{e^{\lambda x_j} + e^{-\lambda x_j}}  
    = \sum_{J\subset\ac{1,\dots,d}} e^{\lambda \pa{\sum_{j \in J} \psca{e_j, x} + \sum_{j \in J^c} \psca{e_j, x}}}
\end{equation*}
we have $\ds e^{\lambda \abs{x}} \le \sum_{J \subset\ac{1,\dots,n}} e^{\lambda \psca{e_J, x}}$ where $(e_J)_j = 1$ if $j \in J$ or $-1$ if $j \in J^c$. With this notation, we have 
\begin{align*}
    \esp{e^{\lambda \abs{\Xmthet}}} &\le \sum_{J\subset\ac{1,\dots,d}} \esp{e^{\lambda \psca{e_J, \Xmthet}}} = \sum_{J\subset\ac{1,\dots,d}} e^{\psi(\lambda e_J-\theta) - \psi(-\theta)}.
\end{align*}
By the  concavity of $\psi-\delta\,|\,.\,|^2$, we have 
\[
\forall\, u,\, v\!\in \R^d,\qquad \psi(u+v)-\psi(u)\le \langle \nabla \psi(u),v\rangle +\delta\,|v|^2
\]
so that
\begin{equation}
    \esp{e^{\lambda \abs{\Xmthet}}} \le \sum_{J\subset\ac{1,\dots,d}} e^{\lambda \psca{\nabla \psi(-\theta), e_J}+\delta\lambda^2|e_J|^2} \le C_{d,\lambda,\delta}  e^{\lambda \sqrt{d} \abs{\nabla \psi(-\theta)}}. \label{esscher-check-NEC-1}
\end{equation}

In the same way, we have 
\begin{align}
    \esp{\abs{\Xmthet}^2 e^{\lambda \abs{\Xmthet}}} &\le \sum_{J\subset\ac{1,\dots,d}} \esp{\abs{X^{(\lambda e_J-\theta)}}^2} e^{\psi(\lambda e_J-\theta) - \psi(-\theta)}, \notag \\ 
    &\le C_{d,\lambda,\delta}  e^{\lambda \sqrt{d} \abs{\nabla \psi(-\theta)}}\sum_{J\subset\ac{1,\dots,d}} \esp{\abs{X^{(\lambda e_J-\theta)}}^2}. \label{esscher-check-NEC-2}
\end{align}
Now, by differentiation of $\psi$ it is easy to check that
\begin{equation*}
    \forall \theta \in \R^d, \quad \Dp \psi(\theta) = \frac{\int x^{\otimes 2} e^{\psca{\theta, x}} p(x) \dd x}{e^{\psi(\theta)}} - \nabla \psi(\theta)^{\otimes 2},
\end{equation*}
which implies 
\begin{equation*}
    \esp{\abs{X^{(\lambda e_J - \theta)}}^2} = \Tr\pa{\Dp \psi(\lambda e_J - \theta)} + \Tr\pa{\nabla \psi(\lambda e_J - \theta)^{\otimes 2}}.
\end{equation*}
The assumption \eqref{Hes} implies that  $0\le D^2\psi(\theta)\le 2\delta\,I_d$ (for the partial order on symmetric matrices induced by nonnegative symmetric matrices) then $D^2\psi(\theta)$ is a bounded function of $\theta\!\in \R^d$ and in turn $\nabla(\theta)$ has a linear growth by the fundamental formula of calculus. Consequently, for every $J \subset \pa{1,\dots,d}$, 
\begin{equation*}
    \esp{\abs{X^{(\lambda e_J - \theta)}}^2} \le C \pa{1 + \abs{\theta}^2}.
\end{equation*}
Plugging this into \eqref{esscher-check-NEC-2} and using \eqref{esscher-check-NEC-1} and \eqref{esscher-check-NEC} we obtain $\esp{\abs{H(\theta, \Xmthet)}^2} \le C\pa{1 + \abs{\theta}^2}$. \hfill $\cqfd$

\section{Adaptive variance reduction for diffusions} \label{diffusion}
\subsection{Framework and preliminaries}
We consider a $d$-dimensional It\^o process $X=(X_t)_{t\in [0,T]}$ solution to the stochastic differential equation (SDE) 
\begin{equation}\label{SDExt} \tag{$E_{b,\sigma,W}$}
    \quad \dd X_t = b(t,X^t) \dd t + \sigma(t,X^t) \dd W_t,\quad X_0=x\!\in \R^d,
\end{equation}
where $W=(W_t)_{t\in [0,T]}$ is a $q$-dimensional standard Brownian motion, $X^t:=(X_{t\wedge s})_{s\in [0,T]}$ is the stopped process at time $t$, $b:[0,T]\times {\cal C}([0,T],\R^d)\to \R^d$ and $\sigma:[0,T]\times  {\cal C}([0,T],\R^d)\to {\cal M}(d,q)$ are measurable with respect to the canonical predictable $\sigma$-field on $[0,T]\times {\cal C}([0,T],\R^d) $. For further details we refer to~\cite{ROWI}, p. 124-130.  

Thus, if $b(t,x^t)=\beta(t,x(t))$ and $\sigma(t,x^t)=\vartheta(t,x(t))$ for every $x\!\in {\cal C}([0,T],\R^d)$, $X$ is a usual diffusion process with drift $\beta$ and diffusion coefficient $\vartheta$.
\ss

If $b(t,x^t)=\beta(t, x(\underline t))$ and $\sigma(t,x^t)=\vartheta(t,x(\underline t))$ for every $x\!\in {\cal C}([0,T],\R^d)$ where $\underline t:= \lfloor \frac {tn}{T}\rfloor \frac Tn $, then $X$ is {\em the continuous Euler scheme with step} $T/n$ of the above diffusion with drift $\beta$ and diffusion coefficient $\vartheta$.

An easy adaptation of standard proofs for regular SDE's show (see~\cite{ROWI}) that strong existence and uniqueness of solutions for~(\ref{SDExt}) follows from the following assumption
\begin{align} \label{Hbs} \tag{${\cal H}_{b, \sigma}$}
    \begin{cases}
    (i)  & b(.,0) \mbox{ and }\sigma(.,0) \mbox { are continuous}\\
    (ii) & \forall\, t\!\in [0,T],\; \forall\, x,\, y\!\in {\cal C}([0,T],\R^d),\, \abs{b(t,y)-b(t,x)}+ \norm{\sigma(t,y)-\sigma(t,x)} \le C_{b,\sigma} \normsup{x-y}.
    \end{cases}
\end{align}

Our aim is to devise an adaptive variance reduction method inspired from Section~\ref{dimfinie} for the computation of
\begin{equation*}
    \esp{F(X)}
\end{equation*}
where $F$ is an Borel functional defined on ${\cal C}([0,T],\R^d)$ such that
\begin{equation}\label{nonvide}
    \prob{F^2(X)>0}>0 \quad \text{and}\quad F(X)\!\in L^2(\P).
\end{equation}
In this functional setting, Girsanov Theorem will play the role of the invariance of Lebesgue measure by translation. The translation process that we consider in this section is of the form $\Theta(t, X^t)$ where $\Theta$ is defined for every $\xi\!\in {\cal C}([0,T],\R^d)$ and $\theta\!\in L^2_{T,p}$  by   
\begin{equation*}
   \Theta(t, \xi) := \varphi(t, \xi^t)\, \theta_t, \quad \text{ with } \quad \varphi : [0,T] \times {\cal C}([0,T],\R^d) \rightarrow {\cal M}(q, p),
\end{equation*}
a bounded Borel function and $\theta\!\in L^2_{T,p}$ (represented by a Borel function) for $p \ge 1$. 
In the sequel, we use the following notations $\varphi_t(\xi) := \varphi(t, \xi^t)$,  
\begin{equation*}
    \Theta_t := \Theta(t, X^t), \quad \Thetath_t := \Theta(t, \Xthett), \quad \text{and} \quad \Thetamth_t := \Theta(t, \Xmthett),
\end{equation*}
where $X^{(\pm \theta)}$ denotes {\color{black} the solution} to $(E_{b \pm \sigma\Theta, \sigma,W})$.

First we need the following standard abstract lemma.
\begin{Lem} \label{Girs} Suppose $(H_{b,\sigma})$ holds. \\
The SDE $(E_{b+\sigma \Theta, \sigma, W})$ satisfies the {\color{black}weak existence and uniqueness} assumptions and for every non negative Borel functional $G:{\cal C}([0,T],\R^{d+1}) \to \R_+$ and $f \in {\cal C}([0,T],\R^{q})$ we have, with the above notations, 
\begin{multline*}
    \esp[4]{G \pa[3]{X, \int_0^.\!\! \psca{f(s, X^s), \dd W_s}}}
    = \esp[4]{G \pa[3]{\Xthet, \int_0^.\!\! \psca{f(s, \Xthets), \dd W_s} 
    + \int_0^. \!\!\psca[1]{f, \Theta}(s, \Xthets) \dd s} \\
    \times e^{- \int_0^T\!\! \psca{\Thetath_s, \dd W_s} -\frac{1}{2} \normLTp{q}{\Thetath}^2}},
\end{multline*}
and
\begin{multline*}
    \esp[4]{G \pa[3]{\Xthet, \int_0^. \!\!\psca{f(s, \Xthets), \dd W_s}}}
    = \esp[4]{G \pa[3]{X, \int_0^. \!\!\psca{f(s, X^s), \dd W_s} 
    - \int_0^. \!\!\psca[1]{f, \Theta}(s, X^s) \dd s}\\
    \times e^{\int_0^T \!\!\psca{\Theta_s, \dd W_s} -\frac{1}{2} \normLTp{q}{\Theta}^2}},
\end{multline*}
\end{Lem}

\ni {\bf Proof.} This is a straightforward application of Theorem 1.11, p.372 (and the remark that immediately follows) in~\cite{REYO} once noticed that $(t,\omega)\mapsto b(t,X^t(\omega))$, $(t,\omega)\mapsto \sigma(t,X^t(\omega))$ and $(t,\omega)\mapsto \Theta(t,X^t(\omega))$ are predictable processes with respect to the completed filtration of $W$.  $\cqfd$


\ni {\bf Remarks.} $\bullet$ The Dol\'eans exponential $\left(e^{\int_0^t\! \psca{\Theta_s, \dd W_s} -\frac{1}{2} \normLTp{q}{\Theta}^2} \right)_{t\in[0,T]}$ is a true martingale for any $\theta\!\in L^2_{T,p}$.

\ni $\bullet$ In fact, still following the above cited remark form~\cite{REYO},  the above lemma holds true if we replace $\Theta$ by any progressively measurable process $\tilde \Theta$ such that $\esp{e^{\frac 12 \int_0^T \!\abs{\tilde \Theta(s,\omega)}^2 \dd s}} < +\infty$.

It follows from the first identity in Lemma~\ref{Girs} that for every bounded Borel function $\varphi:[0,T] \times {\cal C}([0,T], \R^d) \to {\cal M}(q, p)$ and for every $\theta\!\in L^2_{_{T,p}}$
\begin{equation*}
    \esp{F(X)} = \esp{F(\Xthet) 
    e^{- \int_0^T \!\psca{\Thetath_s, \dd W_s} -\frac{1}{2} \normLTp{q}{\Thetath}^2}},
\end{equation*}
(set $G(x,y)=F(x)$). So, finding the estimator with the lowest variance amounts to solving the minimization problem
\begin{equation*}
    \min_{\theta \in L^2_{T,p}} V(\theta) \quad 
    \text{where} \quad V(\theta) := \esp[4]{F^2(\Xthet) e^{-2\int_0^T \!\psca{\Thetath_s, \dd W_s} - \normLTp{q}{\Thetath}^2}},
\end{equation*}
Using Lemma~\ref{Girs} with $G(x,y)=F^2(x) e^{-2 y(T) - \normLTp{q}{\varphi(., x^.) \theta}^2}$ and $f = \Theta$ yields
\begin{equation}\label{V3}
    V(\theta) = \esp[4]{F^2(X) e^{-\int_0^T \!\psca{\Theta_s, \dd W_s} + \frac{1}{2}\normLTp{q}{\Theta}^2}}.
\end{equation}

\begin{Pro} Assume $\E\, F(X)^{2+\eta}<+\infty$ for some $\eta>0$ as well as Assumptions~(\ref{nonvide}) and $(H_{b,\sigma})$. Then function $V$ is finite on $L^2_{T,p}$ and $\log$-convex.
\begin{itemize}
    \item[$(a)$]  Assume that  the bounded matrix-valued Borel function $\varphi$ satisfies that $\varphi(s, X^s)$ has a non-atomic kernel on the event $\{F(X)>0\}$ $i.e.$
\begin{equation}\label{nonatomic}
    \prob{\{\exists \, \theta\!\in L^2_{T,p}\setminus\{0\}\;\mbox{ s.t. } \theta(s) \!\in \mbox{Ker} \, \varphi(s, X^s)\; \dd s\mbox{-a.e. and } F^2(X)>0 \}}=0
    \end{equation}
 then for every finite dimensional subspace $E\subset L^2_{_T}$, $\ds \lim_{\normLTp{p}{\theta} \to +\infty, \theta \in E} V(\theta)=+\infty$.
If furthermore
\begin{equation}\label{General}
    \inf_{\normLTp{p}{\theta} = 1} \int_0^T \theta(s)^*\esp{ \varphi(s, X^s)^*\varphi(s, X^s) \indic[1]{F^2(X) > 0}} \theta(s) \dd s > 0,
\end{equation} 
then $\ds \lim_{\normLTp{p}{\theta} \to +\infty} V(\theta)=+\infty$.
    \item[$(b)$] The function  $V$ is differentiable at every $\theta \!\in L^2_{T,p}$ and the differential $\DD V(\theta) \!\in L^2_{T,p}$  is characterized on every $\psi\!\in L^2_{T,p}$ by
\begin{align} 
    & \pscaLTp{p}{\DD V(\theta), \psi}
    = \esp[4]{F^2(X) e^{-\int_0^T\!\! \psca{\Theta_s, \dd W_s} + \frac{1}{2}\normLTp{q}{\Theta}^2}
    \pa[3]{ \pscaLTp[1]{p}{\Theta, \varphi(., X^.) \psi} - \int_0^T\!\! \psca{\varphi(s, X^s) \psi_s, \dd W_s}}},  \notag \\ 
    & \; = \esp[4]{F^2(\Xmthet) e^{\normLp{L^2_{T,p}}{\Thetamth}^2}
    \pa[3]{2 \pscaLTp[1]{p}{\Thetamth, \varphi(\Xmthet{}^{,.}) \psi} - \int_0^T\!\!  \psca{\varphi(\Xmthets) \psi_s, \dd W_s}}}.\label{Diff2}
\end{align}

\end{itemize}
\end{Pro}

\noindent{\bf Remarks.} $\bullet$ For practical implementation, the ``finite dimensional" statement is the only result of interest since it ensures that $\argmin_{|E}\neq \emptyset$. 

\medskip
\noindent $\bullet$ If $p=q$ and $\varphi=I_q$, the ``infinite-dimensional" assumption is always satisfied.

\bigskip
\ni {\bf Proof.}  $(a)$ As concerns the function $V$, we rely on Equality~(\ref{V3}). Set $r= 1+2/\eta$. Owing to the H\"older Inequality, showing that this function is finite on the whole space $L^2_{T,q}$ amounts to proving that
\begin{equation*} 
    \esp[4]{e^{\frac{r}{2} \normLTp{q}{\Theta}^2 - r\int_0^T \langle\Theta_s ,\dd W_s\rangle}} 
    \le e^{\normsup{\varphi}^2 \normLTp{p}{\theta} r(r+1)/2} < +\infty.
\end{equation*}

   To show that $V$ goes to infinity at infinity, one proceeds as follows. Using the trivial equality 
\begin{equation*}
    e^{-\int_0^T\!\psca{\Theta_s, \dd W_s} + \frac{1}{2}\normLTp{q}{\Theta}^2} = \\
    \pa[4]{e^{-\frac{1}{2} \int_0^T\!\psca{\Theta_s, \dd W_s} + \frac{1}{8}\normLTp{q}{\Theta}^2}}^2
    e^{\frac{1}{4} \normLTp{q}{\Theta}^2}
\end{equation*}
and the reverse H\"older inequality with conjugate exponents $(\frac{1}{3}, -\frac{1}{2})$ we obtain 
\begin{align*}
    V(\theta) & \ge \esp[4]{F^{2/3}(X) e^{\frac{1}{12} \normLTp{q}{\Theta}}}^3 \esp[4]{e^{\frac{1}{2} \int_0^T\!\psca{\Theta_s,\dd W_s} - \frac{1}{8}\normLTp{q}{\Theta}^2}}^{-2}, \\
    & \ge \esp[4]{F^{2/3}(X) e^{\frac{1}{12} \normLTp{q}{\Theta}^2}}^3
\end{align*}
by the martingale property of the Dol\'eans exponential. Let $\varepsilon > 0$ such that $\prob{F^2(X)\ge \varepsilon} > 0$. We have then $V(\theta) \ge \varepsilon^{1/3} \esp[3]{\indic[1]{F^2(X) \ge \varepsilon}e^{\frac{1}{12} \normLTp{q}{\Theta}}}^3$, and by the conditional Jensen inequality 
\begin{align*}
    V(\theta) &\ge \varepsilon^{1/3} \esp[4]{\indic[1]{F^2(X) \ge \varepsilon}e^{\frac{1}{12} \espc[2]{\normLTp{q}{\Theta}^2}{F^2(X) \ge \varepsilon}}}^3\\
    &= \esp[4]{\indic[1]{F^2(X) \ge \varepsilon}e^{\frac{1}{12} \prob{F^2(X) \ge \varepsilon} \esp[2]{\normLTp{q}{\Theta}^2 \indic[1]{F^2(X) \ge \varepsilon}}}}^3 .
\end{align*}

Now
\begin{equation*}
    \esp{\normLTp{q}{\Theta}^2 \indic[1]{F^2(X) \ge \varepsilon}} = \int_0^T \theta(s)^* \esp{\varphi_s(X^s)^* \varphi_s(X^s) \indic[1]{F^2(X) \ge \varepsilon}} \theta(s) \dd s \ge 0.
\end{equation*}
The assumption~(\ref{nonatomic}) implies that, for every $\theta\!\in L^2_{T,p}$,  
\begin{equation*}
    \int_0^T \theta(s)^*\,\esp{\varphi_s(X^s)^* \,\varphi_s(X^s) \indic[1]{F^2(X) \ge \varepsilon}} \theta(s) \dd s \ge \int_0^T \theta(s)^*\,\esp{\varphi_s(X^s)^* \,\varphi_s(X^s) \indic[1]{F^2(X) \ge 0}} \theta(s) \dd s > 0,
\end{equation*} 
so that if $\theta$ runs over the compact sphere of a finite dimensional subspace $E$ of $L^2_{T,p}$
\begin{equation*}
    \inf_{\normLTp{p}{\theta} = 1, \theta\in E}\int_0^T \theta(s)^*\, \esp{\varphi_s(X^s)^*\, \varphi_s(X^s) \indic[1]{F^2(X) \ge \varepsilon}} \theta(s) \dd s > 0,
\end{equation*}
so that
\begin{equation*}
    \lim_{\normLTp{p}{\theta} \to \infty, \theta\in E} \esp[3]{\normLTp{q}{\Theta}^2 \indic[1]{F^2(X) \ge \varepsilon}} = +\infty,
\end{equation*}
and one concludes by Fatou's Lemma using that $\prob{F^2(X) \ge \varepsilon} > 0$.
The second claim easily follows from Assumption~(\ref{General}). 

\ms
\ni $(b)$ As a first step, we show  that the random functional $\Phi(\theta):=\frac{1}{2} \normLTp{q}{\Theta} - \int_0^T \psca{\Theta(s), \dd W_s}$ from  $L^2_{T,p}$ into $L^r(\P)$ ($r\!\in [1,\infty)$),  is differentiable. Indeed, it from the below inequality,
\begin{equation}\label{diffL2Lp}
    \forall\, \theta,\, \psi\!\in L^2_{T, p}, \qquad \abs{\Phi(\theta+\psi) - \Phi(\theta) - \pscaLTp{p}{\DD \Phi(\theta), \psi}} \le \normsup{\varphi}^2 \normLTp{p}{\theta} \normLTp{p}{\psi}
\end{equation}
where $\psi\mapsto \pscaLTp{p}{\DD \Phi(\theta), \psi} = \int_0^T\!\psca{\varphi_s(X^s) \theta(s), \varphi_s(X^s) \psi(s)} \dd s - \int_0^T\!\psca{\varphi_s(X^s) \psi(s), \dd W_s}$ is clearly a bounded random functional from $L^2_{T,p}$ into $L^r(\P)$, with an operator norm $\trnorm{\DD \Phi(\theta)}_{L^2_{T,p}, L^r(\P)} \le \normsup{\varphi}^2 \normLTp{p}{\theta} + c_p \normsup{\varphi}$ ($c_p\!\in (0,+\infty)$ (this follows from H\"older and \emph{B.D.G.} inequalities).

Then, we derive that $\theta \mapsto e^{\Phi(\theta)}$ is differentiable form $L^2_{T,p}$ into every $L^r(\P)$  with differential $e^{\Phi(\theta)} \DD \Phi(\theta)$.  This follows from standard computation based on~(\ref{diffL2Lp}), the elementary inequality $|e^u-1-u|\le \frac 12 u^2(e^{u}+e^{-u})$ and the fact that
\begin{align*}
    \normLp[4]{r}{\abs[3]{\int_0^T \! \psca{\varphi_s(X^s) \psi(s), \dd W_s}} e^{\int_0^T\!\psca{\phi(X^s) \theta(s),\dd W_s}}} 
    & \le \normLp[4]{2r}{\int_0^T\!\psca{\varphi_s(X^s)\psi(s), \dd W_s}} \normLp[4]{2p}{e^{\int_0^T\!\psca{\varphi_s(X^s) \theta(s),\dd W_s}}}, \\
    & \le c_p \normsup{\varphi} \normLTp{p}{\theta} \normLTp{p}{\psi},
\end{align*}
where we used both H\"older and \emph{B.D.G.} inequality.

One concludes that $\theta\mapsto V(\theta)= \esp{F(X)^2e^{\Phi(\theta)}}$ is differentiable by using the $(L^2_{_{T,p}},L^r(\P))$--differentiability of $e^{\Phi(\theta)}$ with $r=1+\frac{\eta}{2}$.  

The second form of the gradient is obtained by a Girsanov transform using Lemma~\ref{Girs}.
$\cqfd$

\subsection{Design of the algorithm}
In view of  a practical implementation of the procedure we are lead to consider some non trivial finite dimensional subspaces $E$ of $L^2_{T,p}$. The function $V$ being strictly $\log$-convex on $E$ and going to infinity as $\normLTp{p}{\theta}$ goes to infinity, $\theta \!\in E$, the restriction of $V$ on $E$ attains a minimum $\theta^*_E$ which {\em de facto} becomes the target of the procedure.  Furthermore, for every $\theta \!\in E$, $\DD V_{|E}(\theta)= \DD V(\theta)_{|E}$ and the quadratic function $L(\theta):= \normLTp{p}{\theta-\theta^*_{_E}}$ is a Lyapunov function for the problem.

Like for the static framework investigated in Section~\ref{translation}, our algorithm will be based on the representation~(\ref{Diff2}) for the differential $\DD V$ of $V$: in this representation the variance reducer $\theta$ appears inside the functional $F$ which makes easier a control at infinity in order to prevent from any early explosion of the procedure. However, to this end we need  to control the discrepancy between $X$ and $\Xmthet$. This is the purpose  of  the following Lemma.

\begin{Lem}  \label{XmoinsXtheta}
Assume $(H_{b,\sigma})$ holds. Let $\varphi$ be a bounded Borel ${\cal M}(q, p)$-valued function defined on $[0,T] \times {\cal C}([0,T], \R^d)$, let $\theta \!\in L^2_{T,p}$ and let $X$ and $X^{(\theta)}$ denote a strong solutions of $E_{b,\sigma,W}$ and $E_{b+\sigma\Theta,\sigma, W}$ driven by the same Brownian motion. Then, for every $r \ge 1$, there exists a real constant $C_{b,\sigma}>0$ such that
\begin{equation} \label{thetasigma}
    \normLp[4]{r}{\sup_{t \in [0,T]} \abs[2]{X_t - \Xthet_t}} \le C_{b,\sigma} e^{C_{b, \sigma} T} \normLp[4]{r}{\int_0^T \!\! \abs[2]{\sigma(s, \Xthets) \Thetath_s} \dd s}.
\end{equation}
\end{Lem}

\ni{\bf Proof.}  The proof follows the lines of the proof of the strong rate of convergence of the Euler scheme (see~$e.g.$~\cite{BOLE}).  $\cqfd$

The main result of this section is the following theorem.
\begin{Thm} Suppose that Assumption~(\ref{nonvide}) and $(H_{b,\sigma})$ hold. 

Let $\varphi$ be a bounded Borel ${\cal M}(q, p)$-valued function (with $p \ge 1$) defined on $[0,T] \times {\cal C}([0,T], \R^d)$, and let $F$ be a functional $F$ satisfying
\begin{equation} \tag{$G_{F,\lambda}$}
    \forall\, x\!\in {\cal C}([0,T],\R^d),\qquad \abs{F(x)} \le C_{_F}(1+ \normsup{x}^\lambda)
\end{equation}
for some positive exponent $\lambda>0$ (then $F(X)\!\in L^r(\P)$ for every $r>0$). Let $E$ be a finite dimensional subspace of $L^2_{T,p}$ spanned by an orthonormal basis $(e_1,\ldots,e_m)$.

Let $\eta > 0$. We define the algorithm by
\begin{equation*}
    \theta_{n+1} = \theta_n -\g_{n+1} H_{\lambda,\eta}(\theta_n,X^{(-\theta_n)}, W^{(n+1)})
\end{equation*}
where $\gamma=(\gamma_n)_{n\ge 1}$ satisfies~(\ref{StepCond}), $(W^{(n)})_{n\ge 1}$ is a sequence of {\color{black}independent Brownian motions for which $X^{(-\theta_n)}={\cal G}(-\theta_n, W^{(n+1)})$ is a strong solution to $(E_{b - \sigma \Theta}, W^{(n+1)})$} and for every standard Brownian motion $W$, every ${\cal F}^W_t$-adapted $\R^p$-valued process $\xi=(\xi_t)_{t\in [0,T]}$,
\begin{equation*}
    \pscaLTp{p}{H_{\lambda,\eta}(\theta,\xi,W), e_i} = 
    \Psi_{\lambda, \eta}(\theta, \xi) F^2(\xi) e^{\normLTp{q}{\Theta(.,\xi^.)}} \pa[3]{2 \pscaLTp[1]{q}{\Theta(.,\xi^.), \varphi(., \xi^.) e_i} - \int_0^T\!\!\psca{\varphi(s, \xi^s) e_i(s),\dd W_s}}
\end{equation*}
where for $\eta > 0$ 
\begin{equation*}
    \Psi_{\lambda, \eta}(\theta, \xi) = \begin{cases}
    \frac{e^{- \normsup{\varphi}\normLTp{p}{\theta}}}{1+\normLTp{q}{\varphi(., \xi^.) \theta}^{2\lambda+\eta}}
    & \text{if $\sigma$ is bounded}, \\
     e^{- (\normsup{\varphi} + \eta)\normLTp{p}{\theta}}    & \text{if $\sigma$ is unbounded}. \\
\end{cases}
\end{equation*}
Then the recursive sequence $(\theta_n)_{n \ge 1}$ a.s. converges toward an $\argmin V$-valued (squared integrable) random variable $\theta^*$.
\end{Thm}
{\color{black}
\ni {\bf Remark.} For a practical implementation of this algorithm, we must have \emph{for all} Brownian motions $W^{(n+1)}$ a strong solution $X^{(-\theta_n)}$ of $(E_{b - \sigma \Theta}, W^{(n+1)})$. In particular, this is the case if the driver $\varphi$ is locally Lipshitz (in space) or if $X$ is the continuous Euler scheme of a diffusion with step $T/n$ (using the driver $\varphi(t, x^t) = f(t, x(\underline t))$).

Note that if $\varphi$ is continuous (in space) but not necessarily locally Lipshitz, the Euler scheme converges in law to the solution of the SDE.
}

\ms
\ni {\bf Proof.}  When the diffusion coefficient $\sigma$ is bounded, it follows from Lemma~\ref{XmoinsXtheta} that, for every $r\ge 1$,
\begin{equation*}
    \normLp[4]{r}{\sup_{t \in [0,T]} \abs[2]{X_t - \Xthet_t}} \le C_{b,\sigma,T} \normsup{\varphi} \normLTp{p}{\theta} \normsup{\sigma},
\end{equation*}
where $\normsup{\sigma} = \sup_{(t,x)\in [0,T]\times{\cal C}([0,T], \R^d)} \norm{\sigma(t,x)}$.

First note that for every $\theta,\,\psi\!\in E$, the mean function $h$ of the algorithm reads
\begin{equation*}
    \pscaLTp{p}{h(\theta), \psi} = \esp{\pscaLTp{p}{H_{\lambda,\eta}(\theta, \Xmthet, W), \psi}} = 
    \esp{\frac{e^{- \normsup{\varphi}\normLTp{p}{\theta}}}{1+\normLTp{q}{\Thetamth}^{2\lambda+\eta}} \pscaLTp{p}{\DD V_{|E}(\theta), \psi}},
\end{equation*}
so that, for every $\theta\neq \theta^*_{_E}$,
\begin{equation*}
    \psca{h(\theta), \theta - \theta^*_{E}} = \esp{\frac{e^{- \normsup{\varphi}\normLTp{p}{\theta}}}{1+\normLTp{q}{\Thetamth}^{2\lambda+\eta}} \pscaLTp{p}{\DD V_{|E}(\theta), \theta-\theta^*_E}} > 0.
\end{equation*}

It remains to check that for every $i\!\in\{1,\ldots,m\}$, $\normLp{2}{H_{\lambda,\eta}(\theta, \Xmthet, W)} \le C\pa[2]{1 + \normLTp{p}{\theta}}$ to apply the Robbins-Zygmund Lemma which ensures the $a.s.$ convergence of the procedure (see Section \ref{argmin-target}).  We first deal with the term $F(\Xmthet)^2\int_0^T\!\psca{\varphi_s(\Xmthets) e_i(s), \dd W_s}$.  Let $\eta'=\frac{\eta}{2\lambda}>0$.
\begin{align*}
    \normLp{2}{F(\Xmthet)^2\int_0^T\!\psca{\varphi_s(\Xmthets) e_i(s), \dd W_s}}
    & \le \normLp{2+\eta'}{F(\Xmthet)^2} \normLp{2(1+1/\eta')}{\int_0^T\!\psca{\varphi_s(\Xmthets) e_i(s), \dd W_s}}, \\
    & \le \normLp{2+\eta'}{F(\Xmthet)^2} \normLp{1+1/\eta'}{\int_0^T\!\abs{\varphi_s(\Xmthets) e_i(s)}^2 \dd s}, \\
    & \le \normLp{2+\eta'}{F(\Xmthet)^2} \normsup{\varphi}.
\end{align*}
Now 
\begin{align*}
    \normLp{2(1+\eta')}{F(\Xmthet)^2} 
    & \le C \pa{1 + \normLp[2]{4\lambda(1+\eta')}{\normsup[2]{\Xmthet}}^{2\lambda(1+\eta')}}, \\
    & \le C_{\lambda,b,\sigma,T} \pa{1 + \normLp[2]{4\lambda(1+\eta')}{\normsup[2]{X}}^{2\lambda(1+\eta')}+ \normLTp{p}{\theta}^{2\lambda(1+\eta')} \normsup{\varphi}^{2\lambda(1+\eta')} \normsup{\sigma}^{2\lambda(1+\eta')}}, \\
    & \le C_{\lambda,b,\sigma,\varphi,T} \pa{1+\normLTp{p}{\theta}^{2\lambda+\eta}}.
\end{align*}
One shows likewise that
\begin{equation*}
    \normLp{2}{F(\Xmthet)^2} \le C_{\lambda,b,\sigma,\varphi, T} \pa{1+\normLTp{p}{\theta}^{2\lambda}}.
\end{equation*}
Combining theses estimates shows that $H_{\lambda,\eta}(\theta, \Xmthet, W)$ satisfies the linear growth assumption in $L^2(\P)$.

\ms
\ni If $\sigma$ is unbounded it follows from Assumption~\eqref{Hbs} that, for every $(t,x)\!\in [0,T]\times {\cal C}([0,T],\R^d)$,
\begin{equation*}
    \norm{\sigma(t,x)} \le C_{\sigma} \pa{1+\normsup{x}}.
\end{equation*}
Elementary computation based on~\eqref{thetasigma} and Lemma~\ref{Girs} yield
\begin{align*}
	\normLp{r}{\int_0^T \abs{\sigma(s, \Xthets) \Thetath_s}} 
	& \le C_\sigma \normLp{L^1_{T,p}}{\theta} \normsup{\varphi} \pa{1 + \normLp{r}{\normsup{X} e^{-\frac{r}{2} \normLTp{q}{\frac{\Theta}{r}}^2 + \int_0^T\!\psca{\frac{\Theta_s}{r}, \dd W_s}}}}, \\ 
	& \le C_\sigma \normLp{L^1_{T,p}}{\theta} \normsup{\varphi} \pa{1 + e^{\frac{\normsup{\varphi}}{2 r r'} \normLTp{p}{\theta}^2} \normLp{r(1+r')}{\normsup{X}}}, \\
	& \le C_{r,b,\sigma,\varphi} \normLTp{p}{\theta},
\end{align*}
for every $r>0$ (Assumption~\eqref{Hbs} implies that $\normLp{r}{\normsup{X}}<+\infty$ for every $r>0$). Following the same proof to the bounded case, we obtain easily the results with $\Psi_{\lambda, \eta}(\theta, \xi) = \frac{e^{- (\normsup{\varphi}+\eta) \normLTp{p}{\theta}}}{1+\normLTp{q}{\varphi(., \xi^.) \theta}^{2\lambda+\eta}}$. We conclude by noting that $\eta$ is an arbitrary parameter to cancel the denominator. $\cqfd$

\ms
\ni {\bf Remark.} If the functional $F$ is bounded ($\lambda = 0$), we prove in the same way that the algorithm without correction, \emph{i.e.} build with $\Psi_{\lambda, \eta} = 1$, a.s. converges. 

\section{Additional remarks} 
For the sake of simplicity we focus in this section on importance sampling by mean translation in a finite dimensional setting (Section~\ref{translation}) although most of the comments below can also be  applied at least  in the path-dependent diffusions setting. 

\subsection{Purely adaptive approach}
As proved by Arouna (see \cite{ARO}), we can consider a purely adaptive approach to reduce the variance. It consists to perform the Robbins-Monro algorithm simultaneously with the Monte Carlo approximation. More precisely, estimate $\esp{F(X)}$ by 
\begin{equation*}
	S_N = \frac{1}{N} \sum_{k=1}^N F(X_k + \theta_{k-1}) \frac{p(X_k + \theta_{k-1})}{p(X_k)}
\end{equation*}
where $X_k$ is the \emph{same innovation} as that used in the Robbins-Monro procedure $\theta_k = \theta_{k-1} - \gamma_k H(\theta_{k-1}, X_k)$.
This adaptive Monte Carlo procedure satisfies a Central Limit Theorem with the optimal asymptotic variance
\begin{equation*}
	\sqrt{N} \pa{S_N - \esp{F(X)}} \xrightarrow{\cal L} {\cal N}(0, \sigma^2_*), \quad \text{whith} \quad \sigma^2_* = V(\theta^*) - \esp{F(X)}^2.
\end{equation*}

This approach can be extended to the Esscher transform when we use the same innovation $\xi_k$ (see \eqref{ReXEsscher}) for the Monte Carlo procedure (computing $X^{(\theta_{k-1})}_k = g(\theta_{k-1}, \xi_k)$) and the Robbins-Monro algorithm (computing $X^{(-\theta_{k-1})}_k = g(-\theta_{k-1}, \xi_k)$). Likewise in the functional setting we can combine the variance reduction procedure  and the Monte Carlo simulations using the same Brownian motion.

In practice, it is not clear that this adaptive Monte Carlo is better than the naive two stage procedure: performing first Robbins-Monro with a small number of iterations (to get a rough estimate $\theta^*$), then performing the Monte Carlo simulations with this optimized parameter. 

\subsection{Weak rate of convergence: Central Limit Theorem (CLT)}
As concerns the rate of convergence, once again this  a regular stochastic algorithm behaves as described in usual Stochastic Approximation Theory textbooks like~\cite{KUYI},~\cite{BEMEPR},~\cite{DUF}. So, as soon as the optimal  variance reducer set is reduced to a single point $\theta^*$, the procedure  satisfies under quite standard assumptions a $CLT$. We will not enter into technicalities at this stage but only try to emphasize the impact of a renormalization factor $g(\theta)$ like $g(\theta):=e^{-\frac{\lambda}{2}|\theta|^b}$ or $g(\theta):=\frac{1}{1+\widetilde F(-\theta)^2}$ induced by the function $F$ on the ``final" rate of convergence of the algorithm toward $\theta^*$. We will assume that $d=1$ and that $X\stackrel{d}{=} {\cal N}(0;1)$ for the sake of simplicity. One can write
\begin{equation*}
    H(\theta,x)= g(\theta)H_0(\theta,x) \quad\mbox{ where }\quad H_0(\theta,x)=F^2(x-\theta) (2\theta-x)
\end{equation*}
The function $H_0$ corresponds to the case of a bounded function $F$ (then $\lambda=0$).
Under simple integration assumptions, one shows that $V$ is twice differentiable and that
\begin{equation*}
    V''(\theta) = e^{\frac{|\theta|^2}{2}}\esp{F^2(X)e^{-\theta X}\pa{1+(\theta-X)^2}}
\end{equation*}

Consequently the mean functions $h$ and $h_0$ related to $H$ and $H_0$ which read respectively
\[
h(\theta)= g(\theta) e^{-|\theta|^2}  V'(\theta) \quad\mbox{ and 
}\quad h_0(x) = e^{-|\theta|^2}  V'(\theta)
\]
are differentiable at $\theta^*$ and
\[
h'(\theta^*)= g(\theta^*) e^{-|\theta^*|^2}  V'(\theta^*)\quad\mbox{ 
and }\quad h_0'(\theta^*)= e^{-|\theta^*|^2}  V'(\theta^*)
\]

Now, general results about CLT say that if $\g_n =\frac{\alpha}{\beta+n}$, $\alpha,\, \beta>0$ with
\begin{equation}\label{CondTCL}
\alpha > \frac{1}{2h'(\theta^*)}= \frac{1}{2g(\theta^*) h_0'(\theta^*)}
\end{equation}
then
\[
\sqrt{n}(\theta_n-\theta^*)\stackrel{{\cal 
L}_{stably}}{\longrightarrow} {\cal N}(0;\Sigma_{\alpha}^*)
\]
where
\begin{equation}\label{asympVar}
\Sigma_{\alpha}^*={\rm Var}(H(y^*,Z))\frac{\a^{2}}{2\a h'(y^*)-1}.
\end{equation}
The mapping $\alpha\mapsto \Sigma_{\alpha}$ reaches its minimum at $\alpha^* = \frac{1}{h'(\theta^*)}=\frac{1}{g(\theta^*) h_0'(\theta^*)}$ leading to the minimal asymptotic variance
\[
\Sigma^*=\Sigma^*_{\alpha^*}= \frac{{\rm Var}(H(y^*,Z))}{h'(y^*)^2}= 
\frac{\esp{H_0(y^*,Z)^2}}{h_0'(y^*)^2}= \frac{\esp{F^4(X)(\theta^*-X)^2e^{-\theta^*X}}}{\esp{F^2(X)(X^2-\theta^*X+1)}^2}
\]
by homogeneity.

So the optimal rate of convergence of the procedure is not impacted by the use of the normalizing function $g(\theta)$.  However, coming back to condition~(\ref{CondTCL}), we see that this assumption on the coefficient $\alpha$ is   more stringent since $\frac{1}{g(\theta^*)}>1$  (in practice this factor can be rather large). Consequently, given the fact that $g(\theta^*)$ is unknown to the user, this will induce a blind choice of $\alpha$ biased to higher values. With the well-known consequence in practice that if $\alpha$ is too large the ``\emph{CLT} regime'' will take place later than it would with smaller values. One solution
to overcome this contradiction can be to make $\alpha$ depend on $n$ and slowly decrease.

As a conclusion,  the algorithm never explodes (and converges) even for strongly unbounded functions $F$ which is a major asset compared to the  version of the  algorithm based on repeated projections. Nevertheless, the normalizing factor which ensures the non-explosion of the procedure may  impact the   rate of convergence since it has an influence on the tuning of the step sequence (which is always more or less ``blind" since it depends on  the target $\theta^*$. In fact, we did not meet such difficulty in our numerical experiments reported below. 

One classical way to overcome this problem can be to introduce the empirical mean of the algorithm implemented with a slowly decreasing step ``\`a la Rupert \& Poliak" (see $e.g.$~\cite{PEL}): Set $\g_n=\frac{c}{n^{r}}$, $\frac 12<r<1$ and
\[
\bar \theta_{n+1} :=\frac{\theta_0+\cdots+\theta_{n}}{n+1}= \bar \theta_n -\frac{1}{n+1}(\bar \theta_n-\theta_{n}),\quad n\ge 0
\]
where $(\theta_n)_{n\ge 0}$ denotes the regular Robbins-Monro algorithm defined by~(\ref{AlgoRM}) starting at $\theta_0$. Then $(\bar \theta_n)_{n\ge 0}$ converges toward $\theta^*$ and satisfies a CLT with the optimal asymptotic variance~(\ref{asympVar}).  See also a variant based on a gliding window developed in~\cite{LEL}.

\subsection{Extension to more general sets of parameters} \label{extension}
In many applications (see below with the Spark spread options with the NIG distribution) the natural set of parameters $\Theta$ is not $\R^q$ but an open connected subset of $\R^q$. Nevertheless, as illustrated below, our unconstrained  approach still works provided one can proceed a diffeomorphic change of parameter by setting
\[
\theta= T(\tilde \theta), \quad \theta \!\in \Theta
\]
where $T: \R^q \to \Theta$  is a ${\cal C}^1$-diffeomorphism with a bounded differential ($i.e.$ $\sup_{\tilde \theta} |\!\| DT(\tilde \theta)|\!\|<+\infty$). As an illustration, let us consider the case where the  state function $H(\theta,X)$ of the procedure  is designed so that $h(\theta):= \E(H(\theta,X))=\rho(\theta) \nabla V(\theta)$ where $V$ is the objective function to be minimized over $\Theta$ and $\rho $ is a bounded {\em positive} Borel function. Then, one  replaces $H(\theta,X)$ by $\widetilde H(\tilde\theta,X):=DT(\tilde\theta).H(T(\tilde \theta),X)$ and defines recursively a procedure on $\R^q$ by
\[
\tilde \theta_{n+1} = \tilde \theta_{n} -\g_{n+1}\widetilde H(\tilde \theta_n,X_{n+1}).
\]
In order  to establish the $a.s.$ convergence of $\theta_n:=T(\tilde \theta_n)$ to $\argmin V$, one  relies on a   variant of Robbins-Monro algorithm, namely a stochastic gradient approach (see~\cite{DUF,KUYI} for further details): one defines  $U(\tilde\theta)= V(T(\tilde \theta))$ which turns out to be a Lyapunov function for the new algorithm since
\[
\langle \nabla U(\tilde \theta),  \E(DT(\tilde\theta)H(T(\tilde \theta),X))\rangle = \rho(T(\tilde \theta)) |\nabla U(\tilde \theta)|^2>0\quad \mbox{ on }\quad T^{-1}(\{\nabla V\neq0\}).
\]
If $U$ satisfies $\|\widetilde H(\tilde \theta,X)\|_{_2}+|\nabla U(\tilde \theta)
|\le C(1+U(\tilde \theta))^{\frac 12}$ (which is a hidden constraint on the choice of $T$), one shows under the standard  ``decreasing" assumption on the step sequence that $U(\tilde \theta_n)\to U_{\infty}\!\in L^1(\P)$ and $\sum_n\g_{n+1} \rho(\tilde \theta_n)  |\nabla U(\tilde \theta_n)|^2<+\infty$. 
If $\ds \lim_{\theta\to \partial \Theta}V(\theta)=+\infty$ or $\ds \liminf_{\theta\to \partial \Theta}\rho(T(\theta))|\nabla V( \theta)|^2>0$, one easily derives that ${\rm dist}(\theta_n, \{\nabla V=0\})  \to 0$ $a.s.$ as $n\to\infty$.

\section{Numerical illustrations} \label{numerical}
\subsection{Multidimensional setting: the NIG distribution}
First we consider a simple case to compare the two algorithms of Section~\ref{dimfinie}. The quantity to compute is 
\begin{equation*}
	\esp{F(X)} = \int_{\R} F(x) \, p_{\NIG}(x;\alpha,\beta,\delta,\mu) \dd x,
\end{equation*}
where $p_{\NIG}(x; \alpha,\beta,\delta,\mu)$ is the density of $X$ a normal inverse gaussian (NIG) random variable of parameters $(\alpha,\beta,\delta,\mu)$ i.e. $\alpha > 0$, $\abs{\beta} \le \alpha$, $\delta > 0$, $\mu \in \R$, 
\begin{equation*}
	p_{\NIG}(x;\alpha,\beta,\delta,\mu) = \frac{\alpha \delta K_1\pa{\alpha \sqrt{\delta^2 +(x-\mu)^2}}}{\pi\sqrt{\delta^2 +(x-\mu)^2}} e^{\delta \gamma + \beta(x-\mu)}, \\
\end{equation*}
where $K_1$ is a modified Bessel function of the second kind and $\gamma = \sqrt{\alpha^2 - \beta^2}$.

We can summarize the two algorithms presented in section \ref{dimfinie}, more precisely the variance reduction based on translation of the density (see Subsection \ref{translation}) and the one based on the Esscher transform (see Subsection \ref{Esscher}), by the following simplified (no computation of the variance) pseudo-code:   
\begin{center} \begin{minipage}{0.48\textwidth}
\centering{Translation (see \ref{translation})}
\begin{lstlisting}[frame= single]
for n = 0 to M do
    X ~ NIG(alpha, beta, mu, delta)
    theta = theta - 1/(n+1000)*H1(theta, X)
for n = 0 to N do
    X ~ NIG(alpha, beta, mu, delta)
    mean = mean + F(X) * p(X+theta)/p(X)
   
\end{lstlisting}
\end{minipage}\hspace{1em} 
\begin{minipage}{0.48\textwidth}
\centering{Esscher transform (see \ref{Esscher})}
\begin{lstlisting}[frame=single]
for n = 0 to M do
    X ~ NIG(alpha, beta-theta, mu, delta)
    theta = theta - 1/(n+1000)*H2(theta, X)
for n = 0 to N do
    X ~ NIG(alpha, beta+theta, mu, delta)
    mean = mean + F(X) * exp(-theta*X)
mean = mean * exp(psi(theta))
\end{lstlisting}
\end{minipage}\end{center}

\begin{itemize} 
\item[--]\emph{Translation case.}
We consider the function $H_1$ of the Robbins-Monro procedure of the first algorithm defined by 
\begin{equation*}
    H_1(\theta, X) = e^{-2\abs{\theta}} F^2(X) \frac{p'(X-2\theta)}{p(X)} \pa{\frac{p(X-\theta)}{p(X-2\theta)}}^2,
\end{equation*}
where an analytic formulation of the derivative $p'$ is easily obtained using the relation on the modified Bessel function $K_1'(x) = \frac{1}{x} K_1(x) - K_2(x)$. 

The assumption \eqref{H2} is satisfied with $a = 1$, and our results of Subsection \ref{translation} apply.

\item[--] \emph{Esscher transform.}
In the Esscher approach we consider the function $H_2$ defined by 
\begin{equation*}
    H_2(\theta, X) = e^{-\abs{\theta}} F^2(X) \pa{\nabla \psi(\theta) - X}.
\end{equation*}
Note that $\psi$ is not well defined for every $\theta \in \R^d$. Indeed, the cumulant generating function of the NIG distribution is defined by 
\begin{equation*}
    \psi(\theta) = \mu \theta + \delta\pa{\gamma - \sqrt{\alpha^2 - (\beta + \theta)^2}},
\end{equation*}
for every $\theta \in (-\alpha-\beta, \alpha-\beta)$. Moreover, we need $\psi(-\theta)$ to be well defined \emph{i.e.} $\theta \in (-\alpha+\beta, \alpha+\beta)$. To take account of these restrictions, we slightly modify the algorithm parametrization (see Subsection \ref{extension}) $\theta = T(\tilde\theta) := (\beta-\alpha) \frac{\tilde \theta}{\sqrt{1+\tilde\theta^2}}$, and update $\tilde{\theta}\in\R$ in the Robbins-Monro procedure (multiply the function $H_2(T(\tilde\theta), X)$ by the derivative $T'(\tilde\theta) = \frac{\beta-\alpha}{(1+\tilde\theta^2)^{3/2}}$). 
\end{itemize}

The payoff $F$ is a Call option of strike $K$, $F(X) = 50(e^X-K)_+$. The parameters of the NIG random variable $X$ are $\alpha=2$, $\beta=0.2$, $\delta=0.8$ and $\mu=0.04$. The variance reduction obtained for different value of $K$ are summarized in the tabular \ref{tab-trans-esscher}. The number of iterations in the Robbins-Monro variance reduction procedure is $M=100\,000$ and the number of Monte Carlo iterations is $N=1\,000\,000$. Note that for each strike, the prices are computed using the same pseudo-random number generator initialized with the same \emph{seed}.
\begin{table}[ht!] 
    \begin{center}
    \begin{tabular}[t]{l|cc|cr|cr}
    \multirow{2}*{K} & \multirow{2}*{mean} & \multirow{2}*{crude var} & var. ratio. & & var. ratio & \\
    & & & translation & ($\theta$) & Esscher & ($\theta$) \\
    \hline
0.6 & 42.19 & 8538 & 5.885 & (0.791) & 56.484 & (1.322)  \\
0.8 & 34.19 & 8388 & 7.525 & (0.903) & 39.797 & (1.309)  \\
1.0 & 27.66 & 8176 & 9.218 & (0.982) & 32.183 & (1.294)  \\
1.2 & 22.60 & 7930 & 10.068 & (1.017) & 29.232 & (1.280) \\
1.4 & 18.76 & 7677 & 9.956 & (1.026) & 28.496 & (1.268)  \\	
    \end{tabular}
    \caption{Variance reduction for different strikes (one dimensional NIG example).}
    \label{tab-trans-esscher} 
\end{center}
\end{table}

\begin{figure}[ht!]
    \begin{center}
        \includegraphics{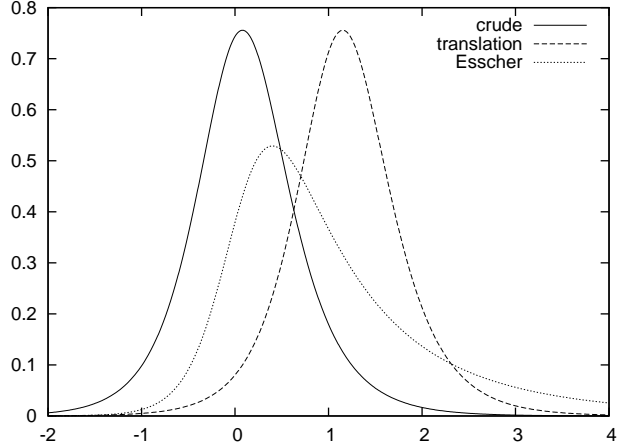}
        \caption{Densities of $X$ (crude), $X+\theta$ (translation) and $\Xthet$ (Esscher) in the case $K=1$.}
        \label{fig-trans-esscher}
    \end{center}
\end{figure}

To complete this numerical example, Figure \ref{fig-trans-esscher} illustrates the densities obtained after the Rob\-bins-Monro procedure. The deformation provided by the Esscher transform is very impressive in this example. We remark that the Esscher transform modifies the parameter $\beta$ which controls the asymmetric shape of the NIG distribution.
\subsubsection*{Spark spread}
We consider now a exchange option between gas and electricity (called spark spread). We choose to model the price of the energy by the exponential of a $NIG$ distribution. 
A simplified form of the payoff is then 
\begin{equation*}
	F(X) = 50 (e^{X^{elec}} - c e^{X^{gas}} - K)_+,
\end{equation*}
where $X^{elec} \sim NIG(2,0.2,0.8,0.04)$ and $X^{gas} \sim NIG(1.4,0.2,0.2,0.04)$ are independent.

The results obtained for different strikes after $300\,000$ iterations of the Robbins-Monro procedure and $3\,000\,000$ iterations of Monte Carlo, are summarized in the Table \ref{tab-spark-spread}. 
\begin{table}[ht!] 
    \begin{center}
    \begin{tabular}[t]{ll|cccc}
    \multirow{2}*{K} & \multirow{2}*{c} & \multirow{2}*{mean} & \multirow{2}*{crude var} & var. ratio. & var. ratio \\
    & & & & translation & Esscher \\
	\hline
	0.4 & 0.2 & 41.021 & 8540.6 & 5.0118 & 25.171 \\
	& 0.4 & 32.719 & 8356.9 & 5.1338 & 27.006 \\
	& 0.6 & 26.337 & 8112.2 & 4.9752 & 28.062 \\
	& 0.8 & 21.556 & 7845.3 & 4.7569 & 29.964 \\
	& 1 & 17.978 & 7582 & 4.5575 & 32.849 \\
	0.6 & 0.2 & 33.235 & 8378.4 & 5.2609 & 27.455 \\
	& 0.4 & 26.534 & 8133.3 & 5.0604 & 28.669 \\
	& 0.6 & 21.587 & 7862.7 & 4.8046 & 30.649 \\
	& 0.8 & 17.931 & 7595.2 & 4.5839 & 33.656 \\
	& 1 & 15.184 & 7344.2 & 4.4064 & 37.489 \\
	0.8 & 0.2 & 26.908 & 8160.1 & 5.1366 & 28.876 \\
	& 0.4 & 21.725 & 7884.9 & 4.844 & 31.018 \\
	& 0.6 & 17.955 & 7612.5 & 4.6031 & 34.166 \\
	& 0.8 & 15.156 & 7357.3 & 4.416 & 38.167 \\
	& 1 & 13.027 & 7123.9 & 4.2685 & 42.781 \\
    \end{tabular}
    \caption{Variance reduction for different strikes (spark spread example).}
    \label{tab-spark-spread} 
\end{center}
\end{table}

\subsection{Functional setting: Down \& In Call option}
We consider a process $(X_t)_{t \ge 0}$ solution of the following diffusion 
\begin{equation*}
    \dd X_t = b(X_t) \dd t + \sigma(X_t) \dd W_t, \quad X_0 = x_0 \in \R.
\end{equation*}
A Down \& In Call option of strike $K$ and barrier $L$ is a Call of strike $K$ which is activated when the underlying $X$ moves down and hits the barrier $L$. The payoff of such a European option is defined by 
\begin{equation*}
    F(X) = (X_T - K)_+ \indic[3]{\ds \min_{0 \le t \le T} X_t \le L}.
\end{equation*}
A naive Monte Carlo approach to price this option is to consider an Euler-Maruyama scheme $\bX=(\bX_\tk)_{k\in\{0,\dots,n\}}$ to discretize $X$ and to approximate $\ds \min_{0 \le t \le T} X_t$ by $\ds \min_{k \in \{0,\dots,n\}} \bX_\tk$. It is well known that this approximation of the functional payoff is poor. More precisely, the weak order of convergence cannot be greater than $\frac{1}{2}$ (see \cite{GOB}). 

A standard approach is to consider the continuous Euler scheme $\bX^c$ obtained by extrapolation of the Brownian between two instants of discretization. More precisely, for every $t \in [\tk, \tkp]$, 
\begin{equation*}
    \bX^c_t = \bX^c_\tk + b(\bX^c_\tk) (t-\tk) + \sigma(\bX^c_\tk)(W_t - W_\tk), \quad \bX^c_0 = x_0 \in \R.
\end{equation*}
By preconditioning, 
\begin{equation} \label{MC_brownian_interpolation}
    \esp{F(X)} = \esp{(\bX_T - K)_+ \pa{1 - \prod_{k = 0}^{N-1} p(\bX_\tk, \bX_\tkp)}},
\end{equation}
with $\ds p(x_k, x_{k+1}) = \probc{\min_{t \in [\tk, \tkp]} \bX^c_t \ge L}{(\bX_\tk, \bX_\tkp) = (x_k, x_{k+1})}$. Now using the Girsanov Theorem and the law of the Brownian bridge (see for example \cite{GLAS}), we have  
\begin{align}
	\begin{split} \label{def-p}
    p(x_k,x_{k+1}) &= 1 - \probc{\min_{t \in [0, t_1]} W_t \le \frac{L-x_k}{\sigma(x_k)}}{W_{t_1} = \frac{x_{k+1}-x_k}{\sigma(x_k)}}, \\
    &= \begin{cases} 0 & \text{if $L\ge\min(x_k, x_{k+1})$}, \\
        1-e^{-\frac{2(L-x_k)(L-x_{k+1})}{\sigma^2(x_k)(\tkp-\tk)}}, & \text{otherwise}.
    \end{cases}
	\end{split}
\end{align}

In the following simulations we consider an Euler scheme of step $\tk = k \frac{T}{n}$ with $n = 100$. 
\subsubsection*{Deterministic case (trivial driver $\varphi \equiv 1$)}
We consider three different basis of $L^2([0,1], \R)$
\begin{itemize}
    \item[--] a polynomial basis composed of the shifted Legendre polynomials $\tilde{P}_n(t)$ defined by
        \begin{equation} \tag{ShLeg}\label{shift_legendre}
    \forall n \ge 0, \forall t \in [0,1], \quad \tilde{P}_n(t) = P_n(2t - 1) \quad \text{where} \quad P_n(t) = \frac{1}{2^n n!} \frac{\dd^{\,n}}{\dd t^n}\pa{(t^2-1)^n}.
\end{equation}
    \item[--] the Karhunen-Lo\`eve basis defined by  
        \begin{equation} \tag{KL} \label{KL}
    \forall n \ge 0, \forall t \in [0,1], \quad e_n(t) = \sqrt{2} \sin \pa{\pa[3]{n+\frac{1}{2}} \pi t} 
\end{equation}
    \item[--] the Haar basis defined by  
        \begin{equation} \tag{Haar} \label{Haar}
    \forall n \ge 0, \forall k = 0,\dots,2^n-1, \forall t \in [0,1], \quad \psi_{n,k}(t) = 2^{\frac{k}{2}} \psi(2^k t - n)
\end{equation}
where
$\psi(t) = \begin{cases} 1 & \text{if $t \in [0, \frac{1}{2})$} \\
        -1 & \text{if $t \in [\frac{1}{2}, 1)$} \\
        0 & \text{otherwise} 
    \end{cases}$
\end{itemize}

\ms
\ni \emph{Black\&Scholes Model} \\
First, we consider the classical Black\&Scholes model. We set the interest rate $r$ to $4\%$ and the volatility $\sigma$ to $70\%$ (which is a high volatility).
The strike of the payoff $F$ is set at $K = 115$ and the barrier level at $L = 65$. A crude Monte Carlo (with Brownian bridge interpolation, see \eqref{MC_brownian_interpolation}) give a price of $2.596$ with a variance of $230$ after $500\,000$ trials. Note that the true price of this product is $2.554$. 

For different basis, the results of our algorithm are summarized in the table \ref{tab-result-bs}. In the Robbins-Monro procedure, we define the step sequence by $\gamma_n = \frac{1}{n+10x_0^2}$ and set the number of iterations at $50\,000$.
\begin{table}[ht!] \label{tab-result-bs}
\centering
\begin{tabular}{lc|ccc}
Basis & Dim. & Mean & ~ CI $95\%$ ~ & Variance ratio \\
\hline
Constant & 1 & 2.5737 & $\pm$0.0230 & 3.4710 \\ 
\hline
ShiftLegendre & 2 & 2.5741 & $\pm$0.0197 & 4.7225 \\ 
\eqref{shift_legendre} & 4 & 2.5717 & $\pm$0.0193 & 4.9478 \\
& 8 & 2.5717 & $\pm$0.0193 & 4.9494 \\
\hline
Karhunen-Lo\`eve & 2 & 2.5678 & $\pm$0.0164 & 6.8644 \\ 
\eqref{KL} & 4 & 2.5729 & $\pm$0.0160 & 7.1851 \\
& 8 & 2.5705 & $\pm$0.0156 & 7.5218 \\
\hline
Haar & 2 & 2.5657 & $\pm$0.0192 & 4.9710 \\ 
\eqref{Haar} & 4 & 2.5671 & $\pm$0.0163 & 6.9459 \\
& 8 & 2.5663 & $\pm$0.0155 & 7.6574 \\
\end{tabular}
\caption{Variance ratio obtained for different basis in the Black\&Scholes model ($K = 115$, $L = 65$, variance of the crude Monte Carlo: 230).}
\end{table}

\begin{figure}[ht!]
    \centering
	 \subfigure[2 vectors]{\includegraphics[width=.49\textwidth]{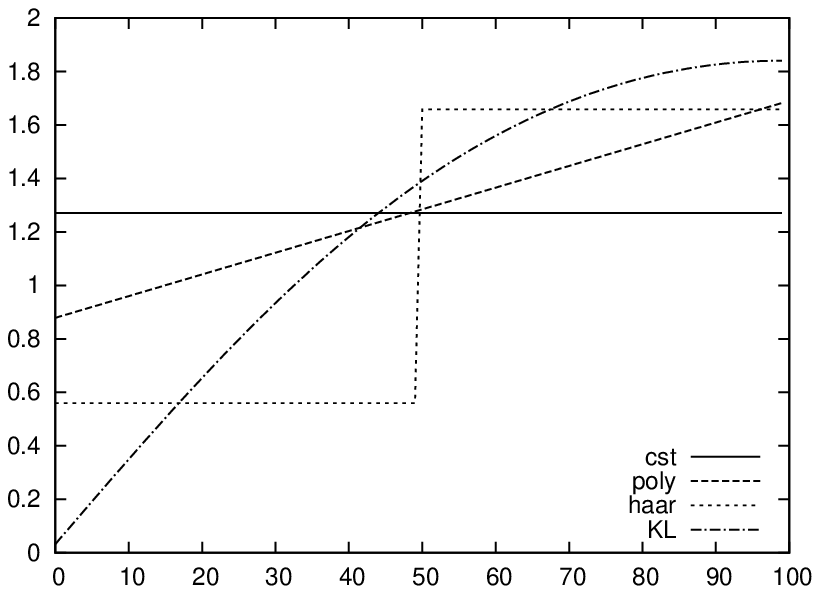}}
    \subfigure[4 vectors]{\includegraphics[width=.49\textwidth]{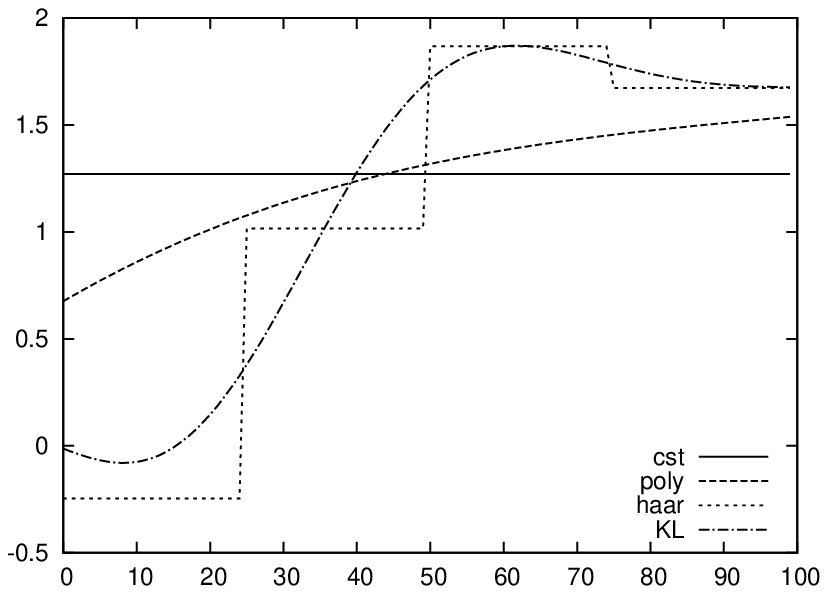}}
    \subfigure[8 vectors]{\includegraphics[width=.49\textwidth]{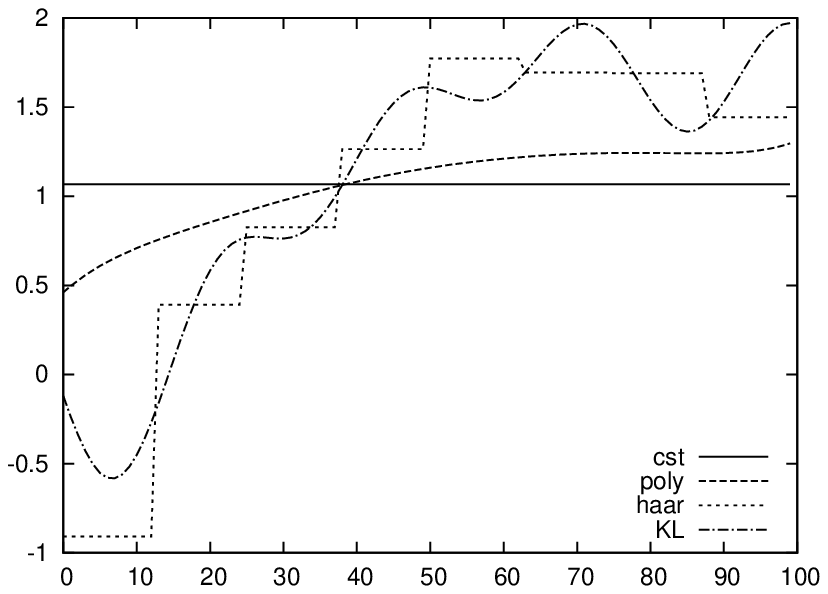}}

\caption{Optimal $\theta$ process obtained with different basis by our algorithm using $50\,000$ trials.} \label{fig-theta}
\end{figure}

In figure \ref{fig-theta} are depicted the optimal variance reducer when the optimization of $V$ is carried out on $E_m$ for several values of $m$ (2, 4 and 8) in the different basis mentioned above.

\ms
\ni \emph{A local volatility Model} \\
To emphasize the generic feature of our algorithm we consider the same product in a local volatility model (inspired by the CEV model) defined by 
\begin{equation} \label{voloc-model}
	\dd x_t = r x_t \dd t + \sigma x_t^\beta \frac{x_t}{\sqrt{1 + x_t^2}} \dd W_t, 
\end{equation}
with $r = 0.04$, $\sigma = 7$ and $\beta = 0.5$.

The price of the Down \& In Call (strike 115, barrier 65) given by a crude Monte Carlo with Brownian interpolation after $500\,000$ trials is $3.194$ and the variance is $206.52$.  
\begin{table}[ht!] \label{tab-result-voloc}
\centering
\begin{tabular}{lc|ccc}
Basis & Dim. & Mean & ~ CI $95\%$ ~& Variance ratio \\
\hline
Constant & 1 & 3.1836 & $\pm$0.0251 & 2.6297 \\
\hline
ShiftLegendre & 2 & 3.1830 & $\pm$0.0223 & 3.3258 \\
\eqref{shift_legendre}& 4 & 3.1815 & $\pm$0.0215 & 3.5670 \\
& 8 & 3.1813 & $\pm$0.0215 & 3.5659 \\
\hline
Karhunen-Lo\`eve & 2 & 3.1852 & $\pm$0.0187 & 4.7254 \\
\eqref{KL} & 4 & 3.1862 & $\pm$0.0183 & 4.9385 \\
& 8 & 3.1918 & $\pm$0.0178 & 5.2183 \\
\hline
Haar & 2 & 3.1834 & $\pm$0.0215 & 3.5699 \\
\eqref{Haar} & 4 & 3.1871 & $\pm$0.0186 & 4.7896 \\
& 8 & 3.1864 & $\pm$0.0177 & 5.2675 \\
\end{tabular}
\caption{Variance ratio obtained for different basis in the local volatility model \eqref{voloc-model} ($K = 115$, $L = 65$, variance of the crude Monte Carlo: 206.52).}
\end{table}

\subsubsection*{Adaptive case (non-trivial driver)}
We experiment now our algorithm with a non-trivial driver  $\varphi$ defined for $t = \tk$ by 
\begin{equation*}
	\varphi(t, \xi^t) = \pa[3]{\bar p_k \quad 1-\bar p_k}, \quad \text{with} \quad \bar p_k = \prod_{j=0}^{k-1} p(\xi_{t_j},\xi_{t_{j+1}}),
\end{equation*}
where $p$ is defined by \eqref{def-p}. Note that $\bar p_k = \probc{\min_{t \in [0, \tk]} \xi_t \ge L}{\xi_0,\dots,\xi_\tk}$ so that there is no extra-computation compared to the Brownian bridge interpolation.

We set $p = 2$ and $E = (\R \ind_{[0,T]})^2$ so that the optimal parameter $\theta_\tk = \alpha \bar p_k + \beta (1-\bar p_k)$ with $(\alpha, \beta) \in \R^2$. The results for different strikes and barrier levels are reported in Table \ref{result-adapt-bs} for the Black\&Scholes model and in Table \ref{result-adapt-voloc} for the local volatility model. The simulation parameters are unchanged.

\begin{table}[ht!] \label{tab-result-adapt-bs}
\centering
\begin{tabular}{cc|ccrr|cc}
Strike & Barrier & Mean & ~ CI 95$\%$ ~ & \multicolumn{2}{c}{Variance ratio (Crude)} & $\alpha$ & $\beta$ \\
\hline
85 & 65 & 2.5738 & $\pm$0.0115 &\hspace{1.5em} 13.49 & (16.56) & -0.1752 & 1.6685 \\
& 75 & 6.0489 & $\pm$0.0186 & 14.26 & (43.39) & 0.0493 & 1.9191 \\
\hline
95 & 65 & 2.5704 & $\pm$0.0110 & 14.64 & (15.26) & 0.0524 & 1.9987 \\
& 75 & 6.0492 & $\pm$0.0190 & 13.67 & (45.25) & 0.1557 & 2.0560 \\
& 85 & 11.5970 & $\pm$0.0301 & 12.23 & (112.92) & 0.4108 & 2.1226 \\
\hline
105 & 65 & 2.5687 & $\pm$0.0122 & 12.03 & (18.56) & 0.3888 & 2.1423 \\
& 75 & 6.0548 & $\pm$0.0206 & 11.66 & (53.08) & 0.3895 & 2.1720 \\
& 85 & 11.5953 & $\pm$0.0308 & 11.67 & (118.32) & 0.4524 & 2.1608 \\
& 95 & 19.2882 & $\pm$0.0348 & 17.17 & (151.04) & 0.6619 & 1.7910 \\
\hline
115 & 65 & 2.5706 & $\pm$0.0135 & 9.75 & (22.90) & 0.5473 & 1.8903 \\
& 75 & 6.0530 & $\pm$0.0211 & 11.16 & (55.42) & 0.4591 & 1.9371 \\
& 85 & 11.5976 & $\pm$0.0297 & 12.55 & (109.98) & 0.4807 & 2.0008 \\
& 95 & 19.2958 & $\pm$0.0347 & 17.21 & (150.67) & 0.7217 & 1.6380 \\
\end{tabular}
\caption{Variance reduction for different strikes and barrier levels in the Black\&Scholes model.} \label{result-adapt-bs}
\end{table}

\begin{table}[ht!] \label{tab-result-adapt-voloc}
\centering
\begin{tabular}{cc|ccrr|cc}
Strike & Barrier & Mean & ~ CI 95$\%$ ~ & \multicolumn{2}{c}{Variance ratio (Crude)} & $\alpha$ & $\beta$ \\
\hline
85 & 65 & 3.1827 & $\pm$0.0127 &\hspace{1.5em}  10.02 & (20.28) & -0.3057 & 1.5522 \\
& 75 & 6.4115 & $\pm$0.0190 & 9.96 & (45.03) & -0.1428 & 1.7985 \\
\hline
95 & 65 & 3.1846 & $\pm$0.0124 & 10.65 & (19.08) & -0.1141 & 1.9139 \\
& 75 & 6.4117 & $\pm$0.0199 & 9.07 & (49.42) & -0.0029 & 1.9814 \\
& 85 & 11.4478 & $\pm$0.0293 & 8.03 & (106.99) & 0.1898 & 1.8937 \\
\hline
105 & 65 & 3.1835 & $\pm$0.0135 & 8.98 & (22.65) & 0.1487 & 1.9628 \\
& 75 & 6.4120 & $\pm$0.0209 & 8.21 & (54.59) & 0.1493 & 2.0060 \\
& 85 & 11.4458 & $\pm$0.0295 & 7.88 & (108.94) & 0.2503 & 1.8737 \\
& 95 & 18.6060 & $\pm$0.0345 & 9.83 & (149.07) & 0.5594 & 1.4343 \\
\hline
115 & 65 & 3.1817 & $\pm$0.0148 & 7.38 & (27.54) & 0.3062 & 1.6884 \\
& 75 & 6.4112 & $\pm$0.0209 & 8.18 & (54.79) & 0.1928 & 1.8119 \\
& 85 & 11.4470 & $\pm$0.0289 & 8.24 & (104.16) & 0.2599 & 1.7430 \\
& 95 & 18.6061 & $\pm$0.0346 & 9.79 & (149.76) & 0.5755 & 1.4313 \\
\end{tabular}
\caption{Variance reduction for different strikes and barrier levels in the local volatility model.} \label{result-adapt-voloc}
\end{table}

\section{Appendix: proof of Theorem~\ref{ThmRZ}}
We propose below the proof of the slight extension of the regular  Robbins-Monro algorithm when $\{h=0\}$ is not reduced to a single equilibrium point.  The key is still   the convergence theorem for non negative super-martingales.

\bs
\ni{\bf Proof.} Set ${\cal F}_n := \sigma(\theta_0,Z_1,\ldots,Z_n)$, $n\ge 1$. Let $\theta^*\!\in {\cal T^*}$. Then
\begin{align}
\notag 
\abs{\theta_{n+1}-\theta^*}^2 &= \abs{\theta_{n}-\theta^*}^2 - 2\g_{n+1} \psca{\theta_{n}-\theta^*, H(\theta_n,Z_{n+1})} + \g_{n+1}^2 \abs{H(\theta_n,Z_{n+1})}^2, \\
\label{IneqL2} & \le \abs{\theta_{n}-\theta^*}^2 - 2\g_{n+1} \psca{\theta_{n}-\theta^*,h(\theta_n)} -2\g_{n+1}\psca{\theta_{n}-\theta^*, \Delta M_{n+1}} + \g_{n+1}^2\abs{H(\theta_n,Z_{n+1})}^2,
\end{align}
where
\begin{equation*}
    \Delta M_{n+1} = H(\theta_n,Z_{n+1}) - \espc{H(\theta_n,Z_{n+1})}{{\cal F}_n} = H(\theta_n,Z_{n+1}) - h(\theta_n),
\end{equation*}
is an increment of (local) martingale satisfying $\esp{|\Delta M_{n+1}|^2} \le C(1+\esp{|\theta_n-\theta^*|^2})$ owing to the assumptions on $H$ and Schwarz Inequality which also implies that
\begin{equation*}
    \esp{\abs{\psca{\theta_{n}-\theta^*, H(\theta_n,Z_{n+1}}}} \le \frac 12 \pa{\esp{\abs{\theta_n-\theta^*}^2} + \esp{\abs{H(\theta_n,Z_{n+1})}^2}} 
    \le C(1+\esp{\abs{\theta_n-\theta^*}^2},
\end{equation*}
for an appropriate real constant $C$. Then, one shows by induction on $n$ from~(\ref{IneqL2}) that $|\theta_n|$ is square integrable for every $n\ge 0$ and that $\Delta M_{n+1}$ is integrable, hence a true martingale increment. 
Now, one derives from the assumptions~(\ref{StepCond}) and~(\ref{IneqL2}) that  
\begin{equation*}
    S_n = \frac{|\theta_n-\theta^*|^2 + 2\sum_{k=0}^{n-1}\g_{k+1} \psca{\theta_{k}-\theta^*, h(\theta_{k})} + C\sum_{k\ge n+1}\g_k^2}{\prod_{k=1}^n (1+C\g^2_k)},
\end{equation*}
is a (non negative) super-martingale with $S_0=|\theta_0-\theta^*|^2\!\in L^1(\P)$. 
This uses  the mean-reverting assumption~\eqref{RMmeanreverting}. Hence $S_n$ is $\P$-$a.s.$ converging toward an integrable r.v. $S_{_\infty}$. Consequently, using that $\sum_{k\ge n+1}\g_k^2\to 0$, one gets
\begin{equation}\label{CvRZ}
    |\theta_n-\theta^*|^2 + 2\sum_{k=0}^{n-1}\g_{k+1} \psca{\theta_k-\theta^*,h(\theta_{k})} \stackrel{a.s.}{\longrightarrow} \widetilde S_{_\infty}=S_{_\infty}\prod_{n\ge1} (1+C_{_L}\g^2_n)\!\in L^1(\P).
\end{equation}
The super-martingale $(S_n)$ being $L^1$-bounded, one derives likewise that $(|\theta_n-\theta^*|^2)_{n\ge 0}$ is $L^1$-bounded since
\begin{equation*}
    |\theta_n-\theta^*|^2\le \prod_{k=1}^n (1+C_{_L}\g^2_k) S_n, \quad n\ge 0.
\end{equation*}
Now, a series with nonnegative terms which is upper bounded by an ($a.s.$) converging sequence, $a.s.$ converges in $\R_+$ so that
\begin{equation*}
    \sum_{n\ge 0} \g_{n+1} \psca{\theta_n-\theta^*, h(\theta_{n})} < +\infty \qquad \P\mbox{-}a.s.
\end{equation*}
It follows from~(\ref{CvRZ}) that, $\P$-$a.s.$, $|\theta_n-\theta^*|^2\stackrel{n\to \infty}{\longrightarrow} L_{_\infty}$ which is integrable since $(|\theta_n-\theta^*|^2)_{n\ge0}$ is $L^1$-bounded and consequently $a.s.$ finite. 

Let $L>0$. Set
\begin{equation*}
\Omega_{_L}:= \left\{\omega\!\in \Omega,\,\forall\, n\ge 0, |\theta_n(\omega)-\theta^*|\le L\right\}.
\end{equation*}
It follows from the $a.s.$ finiteness of $L_{_\infty}$ that $\bigcup_{L>0}\Omega_{_L} = \Omega$ $a.s.$. Now we consider the compact set  $K_{_L}= {\cal T}^*\cap \bar B(0,L)$. It is separable so there exists an everywhere dense  sequence in $K_{_L}$, denoted for convenience $(\theta^{*,k})_{k\ge 1}$. The above proof shows that $\P$-$a.s.$, for every $k\ge 1$,  $|\theta_n-\theta^{*,k}|^2\to \ L^k_{_\infty}<+\infty$ as $n\to \infty$. Then set
\begin{equation*}
    \Omega'_{_L}:=\left\{\omega\!\in \Omega_{_L},\; \abs{\theta_n(\omega)-\theta^{*,k}}^2 \stackrel{n\to \infty}{\to} L^k_{_\infty}(\omega), k\ge 1, \; \sum_{n\ge 1}\gamma_n \psca{\theta_{n-1}(\omega)-\theta^*, h(\theta_{n-1}(\omega)} <+ \infty\right\}
\end{equation*}
which satisfies $\P(\Omega'_{_L})= \P(\Omega_{_L})$. Assume $\omega\!\in \Omega'_{_L}$. Up to  two successive extractions, there exists a subsequence $\theta_{\phi(n,\omega)}$ such that 
\begin{equation*}
    \psca{\theta_{\phi(n,\omega)}-\theta^{*}, h(\theta_{\phi(n,\omega)}(\omega))} \stackrel{n\to \infty}{\longrightarrow} 0\qquad\mbox{ and } \qquad \theta_{\phi(n,\omega)}(\omega) \stackrel{n\to \infty}{\longrightarrow}\theta_{_\infty}(\omega).
\end{equation*}
The function $h$ being continuous $\psca{\theta_{_\infty}(\omega)-\theta^{*}, h(\theta_{_\infty}(\omega))}=0$ which implies that $\theta_{_\infty}(\omega)\!\in\{h=0\}$. Hence $\theta_{_\infty}(\omega)\!\in K_{_L}$. Then any limiting value $\theta'_{_\infty}(\omega)$ of the sequence $(\theta_n(\omega))_{n\ge 1}$ will satisfy
\begin{equation*}
    \forall\, k\ge 1,\quad |\theta'_{_\infty}(\omega)-\theta^{*,k}|= |\theta_{_\infty}(\omega)-\theta^{*,k}|= \sqrt{L^k_{_\infty}(\omega)}
\end{equation*}
which in turn implies that $\theta'_{_\infty}(\omega)=\theta_{_\infty}(\omega)$ by considering a subsequence $\theta^{*,k'}\to \theta_{_\infty}(\omega)$. So, $\theta_{_\infty}(\omega)$ is the unique limiting value of the sequence $(\theta_n(\omega))_{n\ge 0}$ $i.e.$ $\theta_n(\omega)\to \theta_{_\infty}(\omega)$ as $n\to \infty$.  The fact that the resulting random vector $\theta_{_\infty}$ is   square integrable follows from  Fatou's Lemma and the $L^2$-boundedness of the sequence $(\theta_n-\theta^*)_{n\ge 1}$.$\cqfd$
\end{document}